\documentclass{article}
\usepackage{amsmath,amsrefs, amssymb, amsthm, mathtools}
\usepackage{afterpage}
\usepackage{accents}
\usepackage[shortlabels]{enumitem}
\usepackage{float}
\usepackage{geometry}
\usepackage{hyperref}
\usepackage{subcaption}
\usepackage{tikz}
\usepackage{xcolor}
\usetikzlibrary{arrows.meta}
\newtheorem{theorem}{Theorem}[section]
\newtheorem{lemma}[theorem]{Lemma}

\newtheorem{conjecture}[theorem]{Conjecture}
\newtheorem{corollary}[theorem]{Corollary}
\newtheorem{proposition}[theorem]{Proposition}

\theoremstyle{definition}
\newtheorem{definition}[theorem]{Definition}
\newtheorem{remark}[theorem]{Remark}
\newtheorem{observations}[theorem]{Observations}
\newtheorem{example}[theorem]{Example}

\newcommand{\N}{\mathbb{N}}

\newcommand{\LSQ}{\mathsf{LSQ}}
\newcommand{\LD}{\mathsf{LD}}
\newcommand{\stLSQ}{\mathsf{stLSQ}}
\newcommand{\stLD}{\mathsf{stLD}}

\newcommand{\CC}{\mathcal{CC}}

\newcommand{\dinv}{\mathsf{dinv}}
\newcommand{\Attack}{\mathsf{Attack}}
\newcommand{\area}{\mathsf{area}}
\newcommand{\shift}{\mathsf{shift}}
\newcommand{\dw}{\mathsf{dw}}
\newcommand{\sdw}{\mathsf{sdw}}
\newcommand{\sched}{\mathsf{sched}}
\newcommand{\maj}{\mathsf{maj}}
\newcommand{\revmaj}{\mathsf{revmaj}}
\newcommand{\ADR}{\mathsf{ADR}}
\newcommand{\DADR}{\mathsf{DADR}}
\newcommand{\paritydec}{\mathsf{parity-dec}}

\newcommand*{\dec}[1]{%
    \accentset{\mbox{\LARGE\bfseries .}}{#1}
    }

\title{Decorated square paths at $q=-1$}
\author{Sylvie Corteel and Alexander Lazar and Anna Vanden Wyngaerd}
\begin{document}
\maketitle
\begin{abstract}
    The \emph{valley delta square} conjecture states that the symmetric function $\frac{[n-k]_q}{[n]_q}\Delta_{e_{n-k}}\omega(p_n)$ can be expressed as the enumerator of a certain class of \emph{decorated square paths} with respect to the bistatistic $(\dinv,\area)$. Inspired by recent positivity results of Corteel, Josuat-Verg\`{e}s, and Vanden Wyngaerd, we study the evaluation of this enumerator at $q=-1$. By considering a cyclic group action on the decorated square paths which we call \emph{cutting and pasting}, we show that $\left.\left\langle \frac{[n-k]_q}{[n]_q}\Delta_{e_{n-k}}\omega(p_n), h_1^n\right\rangle\right|_{q=-1}$ is $0$ whenever $n-k$ is even, and is a positive polynomial related to the Euler numbers when $n-k$ is odd. We also show that the combinatorics of this enumerator is closely connected to that of the Dyck path enumerator for $\langle\Delta_{e_{n-k-1}}'e_n,h_1^n\rangle$ considered by Corteel--Josuat Verg\`{e}s--Vanden Wyngaerd.
\end{abstract}
\tableofcontents
\section{Introduction}

The symmetric function $\nabla e_n$, also known as the Frobenius characteristic of the ring of diagonal coinvariants \cite{Haiman2002}, has been the subject of a multitude of papers in algebraic combinatorics and related fields in the last three decades. In \cite{HaglundHaimanLoehrRemmelUlyanov2005}, the authors proposed a combinatorial formula of this remarkable function in terms of \emph{decorated Dyck paths}, which was proved a decade later in \cite{CarlssonMellit2018}. This result is known as the \emph{shuffle theorem}. Many special cases of this formula lead to interesting combinatorics: $\langle \nabla e_n,e_n\rangle$, $\langle \nabla e_n,e_{n-d}h_d\rangle$, $\langle \nabla e_n,h_{1^{n}}\rangle$ are the $q,t$-Catalan numbers, $q,t$-Schröder numbers and $q,t$-parking functions, respectively (see \cites{GarsiaHaglund2002,Haglund2004}). 

Many generalizations and analogues to the shuffle formula have been studied. In \cite{HaglundRemmelWilson2018}, the authors provide conjectural formulas for the symmetric function $\Delta'_{e_{n-k-1}}e_n$, in terms of \emph{decorated Dyck paths}, which reduces to the shuffle theorem when $k=0$.   
In this paper, we will in particular be concerned with the Hilbert series of the \emph{valley version} of the \emph{Delta conjecture}, which can be stated as follows
$$\langle\Delta'_{e_{n-k-1}} e_n,h_1^n\rangle = \sum_{P \in \stLD(n)^{\bullet k}}t^{\area(P)}q^{\dinv(P)},$$
where the sum is over \emph{valley-decorated standardly labeled Dyck paths} (see Section~\ref{sec:valley-delta-square} for all the precise definitions). 

Another result related to this story is the \emph{square theorem}, which gives a formula for $\nabla \omega p_n$ in terms of \emph{square paths} (conjectured in \cite{LoehrWarrington2007} and proved in \cite{Sergel2016}). A ``Delta generalization'' of the square theorem was proposed in \cite{DAdderioIraciVandenWyngaerd2020deltasquare}. Of interest to us is again its Hilbert series: 
\[\frac{[n-k]_q}{[n]_q}\left\langle\Delta_{e_{n-k}} \omega(p_n),h_{1}^n \right\rangle= \sum_{P\in \stLSQ(n)^{\bullet k}}q^{\dinv(P)}t^{\area(P)},\] where the sum is over \emph{valley-decorated standardly labeled square paths} (we refer again to Section~\ref{sec:valley-delta-square} for precise definitions). 

As these formulas have become more intricate, a number of elegant results have been found when specializing one of the $q,t$-variables to $0$ or $1$ \cites{IraciRhoadesRomero2022,RhoadesWilson2023,IraciNadeauVandenWyngaerd2023}.

In \cite{CorteelJosuatVergesVandenWyngaerd2023}, the authors study
$$ D_{n,k} \coloneqq \sum_{P \in \stLD(n)^{\bullet K}}(-1)^{\dinv(P)}t^{\area(P)},$$
i.e.\ the specialization of the combinatorial enumerator of the valley Delta conjecture at $q=-1$, which has nicer properties than one might expect.
They prove that
\begin{equation}
    \sum_{k=0}^{n-1}D_{n,k}z^k = \sum_{\sigma\in \mathfrak{S}_n}t^{\mathsf{inv}_3(\sigma)}z^{\mathsf{monot}(\sigma)},\label{eq:sum-Dnkz}
\end{equation}
where $\mathrm{inv}_3$ and $\mathrm{monot}$ are certain permutation statistics. Taking $z=0$, and via the shuffle theorem, this becomes 
$$\left.\langle\nabla e_n,h_{1^n}\rangle\right|_{q=-1}=D_{n,0}= t^{\lfloor n^2/4\rfloor}E_n(t),$$
where $E_n(t)$ is a known $t$-analog of the Euler numbers introduced in \cite{HanRandrianarivonyZeng1999} and further studied in \cites{Chebikin2008, JosuatVerges2010}: it $t$-counts alternating permutations\footnote{A permutation $\sigma$ is said to be \emph{alternating} if $\sigma_1>\sigma_2<\sigma_3>\cdots$. It is well-known that these permutations are counted by the Euler numbers.} with respect to $31-2$ patterns\footnote{A $31-2$-pattern of a permutation $\sigma\in \mathfrak S_n$ is a pair $(i,j)$ such that $1<i+1<j\leq n$ and $\sigma_{i+1}<\sigma_j<\sigma_i$.}. 
This polynomial identity is itself a $t$-refinement of a more classical identity relating the Euler numbers $E_n$ to an alternating sum over the set of parking functions (c.f. the introduction of \cite{CorteelJosuatVergesVandenWyngaerd2023}).

Inspired by \cite{CorteelJosuatVergesVandenWyngaerd2023}, in this work we study 
$$ S_{n,k} \coloneqq \sum_{P \in \stLSQ(n)^{\bullet k}}(-1)^{\dinv(P)}t^{\area(P)},$$
i.e.\ the square paths enumerator in the valley delta square conjecture, and find that it \emph{also} has well-behaved combinatorics at $q=-1$.

We show that $S_{n,k} = 0$ when $n-k$ is even and that 
\begin{equation} 
    \sum_{k=0}^{n-1}S_{n,k}z^{k} = \sum_{\sigma\in\mathfrak S_n} t^{\revmaj(\sigma)}z^{\paritydec(\sigma)},\label{eq:sum-Snkz}
\end{equation} 
where $\paritydec$ is a combinatorial quantity (it is the number of decorations resulting from the parity decorating algorithm, see Proposition~\ref{prop:decorating-algos}).
The $z=0$ case, combined with the square theorem, becomes
$$\left.\langle \nabla \omega p_n, h_{1^n} \rangle\right|_{q=-1} = S_{n,0} = \begin{cases}[n]_t t^{\lfloor (n-1)^2/4\rfloor}E_{n-1}(t), &  \text{if $n$ is odd}, \\ 0, & \text{if $n$ is even}\end{cases},$$
and specializing to $t=1$ yields
$$\sum_{P \in \stLSQ(n)}(-1)^{\dinv(P)} = \begin{cases} nE_{n-1}, & \text{if $n$ is odd}\\ 0, & \text{if $n$ is even}\end{cases},$$
which are the coefficients of $\frac{x^n}{n!}$ in the Taylor series of $\frac{x}{\cos(x)}$ (OEIS sequence A009843). 

The strategy we use for proving \eqref{eq:sum-Snkz} starts out as very similar to the strategy for proving \eqref{eq:sum-Dnkz}: we use the powerful tool of \emph{schedule numbers} (see Section~\ref{sec:schedule}), to reduce the enumeration over all paths to paths that have schedule numbers all $1$. However, in the square paths case, this does not yet yield a positive enumeration. We explain the remaining cancellations by defining an equivalence relation on the schedule $1$ paths of ``being in the same \emph{cutting cycle}''. We then show that each of these equivalence classes contributes either $1$ (when $n-k$ is odd), or $0$ (when $n-k$ is even), to the enumeration in $S_{n,k}$. Finally, we show that there is a one to one correspondence between permutations of $n$ and equivalence classes of square paths of size $n$ with $k$ decorations for which $n-k$ is odd.  

Furthermore, we discuss the strong connection between the valley Delta conjecture and the valley Delta square conjecture at $q=-1$.
On the symmetric function side, by Theorems 4.11 and 4.12 in \cite{DAdderioIraciVandenWyngaerd2020deltasquare}, we know that
\begin{align*}
    \left.\sum_{k=0}^{n-1}\frac{[n-k]_q}{[n]_q}\Delta_{e_{n-k}}\omega(p_n)\right|_{q=-1} = \left.\sum_{k=0}^{n-1} \Delta'_{e_{n-k-1}}e_n\right|_{q=-1} =  \left.\nabla e_n\right|_{q=0}\label{eq:sf-sum-k}
\end{align*}
Thus, taking the scalar product with $h_{\lambda}$ yields the multinomial coefficient $\genfrac[]{0pt}{1}{n}{\lambda}_q.$
When $\lambda = 1^n$, we will establish on the combinatorial side that \[\sum_{k=0}^{n-1} S_{n,k} = \sum_{k=0}^{n-1} D_{n,k} = [n]_t!\] by exhibiting an explicit bijection of the combinatorial enumerators of both sides.
Next, we show that we have the following recursive relationship
$$ S_{n,k} = \begin{cases} [n]_t(D_{n-1,k} + D_{n-1,k-1}), & \text{if $n-k$ is odd}\\ 0, & \text{if $n-k$ is even} \end{cases},$$
where $D_{-1,k} = D_{n,-1} = 0$. 

The rest of the paper is organized as follows: in Section~\ref{sec:valley-delta-square}, we define all the necessary combinatorics of lattice paths to state the Valley Delta square conjecture. In Section~\ref{sec:schedule} we recall the notion of \emph{schedule numbers}, which will be a crucial tool in our proofs. Section~\ref{sec:enumeration} contains the motivation and statement of our enumerative formula for $S_{n,k}$, and its proof is developed in Section~\ref{sec:cutting}. In Section~\ref{sec:dycktosquare} we explore the relationship between the enumerators of decorated Dyck paths and square paths at $q=-1$. 

\section{Valley Delta square conjecture}\label{sec:valley-delta-square}

In this section we give the necessary combinatorial definitions to state the valley Delta square conjecture. These definitions were introduced in \cites{HaglundRemmelWilson2018,DAdderioIraciVandenWyngaerd2020deltasquare}.

\begin{figure}
    \centering
    \begin{tikzpicture}[scale = .5]
        \draw[step=1.0, gray!60, thin] (0,0) grid (7,7);

        \draw[gray!60, thin] (2,0) -- (7,5);

        \draw[blue!60, line width=3pt] (0,0) -- (1,0) -- (1,1) -- (2,1) -- (3,1) -- (3,2) -- (3,3) -- (3,4) -- (4,4) -- (4,5) -- (5,5) -- (5,6) -- (6,6) -- (6,7) -- (7,7);

        \node at (1.5,0.5) {$2$};
        \draw (1.5,0.5) circle (.4cm);
        \node at (3.5,1.5) {$1$};
        \draw (3.5,1.5) circle (.4cm);
        \node at (3.5,2.5) {$2$};
        \draw (3.5,2.5) circle (.4cm);
        \node at (3.5,3.5) {$3$};
        \draw (3.5,3.5) circle (.4cm);
        \node at (4.5,4.5) {$1$};
        \draw (4.5,4.5) circle (.4cm);
        \node at (5.5,5.5) {$2$};
        \draw (5.5,5.5) circle (.4cm);
        \node at (6.5,6.5) {$3$};
        \draw (6.5,6.5) circle (.4cm);
        \node at (4.5,5.5) {$\bullet$};
        \node at (1--2-0.5,1+0.5) {$\bullet$};
        \node at (6-0-0.5,6+0.5) {$\bullet$};
    \end{tikzpicture}
    \caption{An element of $\LSQ(7)^{\bullet 3}$.}\label{fig:squarepath}
\end{figure}
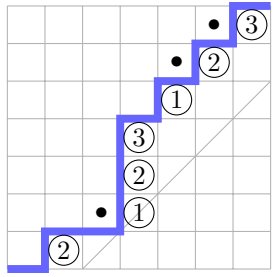

Throughout this paper, $n$ will be a positive integer.
\begin{definition}
    A \emph{square path} of size $n$ is a lattice path in the plane
    \begin{itemize}
        \item consisting of unit north and east steps (also called \emph{vertical} and \emph{horizontal} steps),
        \item starting at $(0,0)$ and ending at  $(n,n)$,
        \item ending with an east step.
    \end{itemize}
    A square path is called a \emph{Dyck path} if it stays weakly above \emph{main diagonal}, that is the line $x=y$.
\end{definition}
Let us fix some useful vocabulary for talking about lattice paths.
\begin{itemize}
    \item We call a \emph{square} any area one square in the plane whose vertices have integer coordinates.
    \item Furthermore, we will always consider steps of square paths to be either left or top edges of squares. This way, we can refer to \emph{the} square of a step without ambiguity.
    \item We refer to the collection of squares that have two vertices on the line $y=x+i$ as the \emph{$i$-th diagonal}.
    \item The $0$-diagonal is sometimes referred to as the \emph{main diagonal}. The \emph{bottom diagonal} of a path is the lowest diagonal that intersects with the path.
\end{itemize}

\begin{definition}
    Given a square path $\pi$, a \emph{labeling} of $\pi$ is a word $w$ of positive integers whose $i$-th letter labels the $i$-th vertical step of $\pi$ such that, when the label is placed in the square of its step, the labels appearing in squares of the same column are increasing from bottom to top.
    Such a pair $(\pi,w)$ is called a \emph{labeled square path}.  The set of labeled square (respectively, Dyck) paths of size $n$ is denoted by $\LSQ(n)$  (respectively $\LD(n)$).
\end{definition}
\begin{definition}
    A labeling of a path of size $n$ is said to be \emph{standard} if its labels are exactly $1,2,\dots, n$. The set of standardly labeled square (respectively, Dyck) paths is denoted $\stLSQ(n)$ (respectively, $\stLD(n)$).
\end{definition}
The reader might be more familiar with the set of standardly labeled Dyck paths of size $n$ as the set of \emph{parking functions} of the same size. 

\begin{definition}
    The \emph{area word} of a square path of size $n$ is the word $a$ of $n$ integers whose $i$-th letter is the index of the diagonal of the $i$-th vertical step, i.e.\ it equals $j$ if the starting point of the $i$-th vertical step lies on the line $y=x+j$.
\end{definition}

\begin{definition}
    The \emph{shift} of a square path $\pi$ of area word $a$, is the absolute value of the minimum letter of its area word. Notice that the shift of a Dyck path is always $0$.
\end{definition}
\begin{definition}
    The \emph{area} of a square path $\pi$ of area word $a$ and shift $s$ is \[\area(\pi) = \sum_i (a_i + s).\] For a labeled square path $(\pi,w)$ we define its area to be independent of its labeling: \[\area((\pi,w)) \coloneqq \area(\pi).\] In other words, the area of a path is the number of whole squares between the paths and the line $y=x-s$.
\end{definition}

\begin{example}
    The path in Figure~\ref{fig:squarepath} has area word $(-1,-2,-1,0,0,0,0)$, its shift is $2$ and its area is $10$.
\end{example}

\begin{definition}
    Given $P\coloneqq (\pi,w)\in \LSQ(n)$ with area word $a$, the $i$-th vertical step of $P$ is called a \emph{contractible valley} if
    \begin{itemize}
        \item either $i>1$ and $a_{i-1}> a_{i}$,
        \item or $i>1$ and $a_{i-1} = a_{i}$ and $w_{i-1} < w_{i}$,
        \item or $i=1$ and $a_i\leq -1$.
    \end{itemize}
    In other words, the $i$-th vertical step is a contractible valley if it is preceded by a horizontal step and the following holds:  after replacing the two steps
    \tikz{\draw[blue!60, line width=1.5pt] (0,0)-|(.3,.3);} with \tikz{\draw[blue!60, line width=1.5pt] (0,0)|-(.3,.3);} (and accordingly shifting the $i$-th label one cell to the left), we still get a valid labeled path where labels are increasing in each column.
\end{definition}

\begin{definition}
    A \emph{(valley) decorated labeled square (respectively, Dyck) path} is a triple $(\pi,w,dv)$ where $(\pi,w)\in \LSQ(n)$ (respectively $\LD(n)$) and $dv$ is some subset of the contractible valleys of $(\pi,w)$. The elements of $dv$ are called \emph{decorations}, and we visualize them by drawing a $\bullet$ to the left of these contractible valleys. The set $\LD(n)^{\bullet k}$ denotes the decorated labeled Dyck paths with exactly $k$ decorations.
\end{definition}

\begin{example}
    The path in Figure~\ref{fig:squarepath} has $4$ contractible valleys: the vertical steps of indices $1,2,6,7$. Three of them ($2,6$ and $7$) have been decorated.
\end{example}

\begin{definition}
    Let $P\coloneqq (\pi,w, dv) \in \LSQ(n)^{\bullet k}$ with area word $a$ and $(i,j)$ a pair of indices of vertical steps with $1 \leq i < j\leq n$. These steps are said to \emph{attack each other}, or to be in an \emph{attack relation} if
    \begin{itemize}
        \item either $a_i = a_j, w_i<w_j$ and  $i\not \in dv$,
        \item or $a_i = a_j+1, w_i>w_j$ and  $i\not \in dv$.
    \end{itemize}
    The set of such pairs of indices is denoted by $\Attack(P)$. An attack relation of the first kind is referred to as \emph{primary dinv} and of the second kind as \emph{secondary dinv}.
\end{definition}

\begin{definition}
    Given a path $P\in \LSQ(n)^{\bullet k}$ with area word $a$, we define its \emph{dinv} to be \[\dinv(P)\coloneqq \#\Attack(P) + \#\{i\mid a_i<0\} - k.\]
    The second term of this sum is the number of labels in negative diagonals and is referred to as \emph{tertiary} or \emph{bonus} dinv.
\end{definition}
\begin{remark}\label{rem:dinv-nonneg}
    Though this is not immediate, it is always the case that $\dinv(P)\geq 0$. Indeed, the presence of a contractible valley always necessarily implies the presence of a unit of primary, secondary or tertiary dinv. See \cite{VandenWyngaerd2020}*{Proposition~4.20}
\end{remark}

\begin{example}
    The attack relations for the path in Figure~\ref{fig:squarepath} are $(1,2),(5,6)$ and $(5,7)$ (the first is secondary and the other two primary dinv). There are $3$ labels under the line $x=y$, so $3$ units of tertiary dinv. The number of decorated valleys is $3$ so the total dinv comes to $3+3-3=3$.
\end{example}

\begin{definition}
    Given $P = (\pi,w,dv)\in\LSQ(n)^{\bullet k}$, we define the monomial \[x^P = x_{w} \coloneqq \prod_i x_{w_i}.\]
\end{definition}
\begin{example}
    The monomial for the path in Figure~\ref{fig:squarepath} is $x_1^2x_2^3x_3^2$.
\end{example}

\begin{conjecture}[The valley Delta conjecture \cite{HaglundRemmelWilson2018}]\label{conj:valley-delta}
    For all $n,k\in \N$, \[\Delta'_{e_{n-k-1}}e_n = \sum_{P\in \LD(n)^{\bullet k}}q^{\dinv(P)}t^{\area(P)}x^P\]
\end{conjecture}
\begin{conjecture}[The valley Delta square conjecture \cite{IraciVandenWyngaerd2021}]\label{conj:valley-delta-square}
    For all $n,k\in\N$,
    \[\frac{[n-k]_q}{[n]_q}\Delta_{e_{n-k}} \omega(p_n) = \sum_{P\in \LSQ(n)^{\bullet k}}q^{\dinv(P)}t^{\area(P)}x^P.\]
\end{conjecture}

We will in particular be interested in the Hilbert series of these expressions. Taking the scalar product $\langle\cdot,h_1^{n}\rangle$:

\begin{align}
    \langle \Delta'_{e_{n-k-1}}e_n,h_1^{n}\rangle                            & = \sum_{P\in \stLD(n)^{\bullet k}}q^{\dinv(P)}t^{\area(P)}\label{eq:standard-delta}          \\
    \frac{[n-k]_q}{[n]_q}\langle \Delta_{e_{n-k}} \omega(p_n),h_1^{n}\rangle & = \sum_{P\in \stLSQ(n)^{\bullet k}}q^{\dinv(P)}t^{\area(P)}\label{eq:standard-delta-square}.
\end{align}

In this paper, we will study the combinatorial side of these equations, evaluated at $q=-1$:
\begin{align}
    D_{n,k} & \coloneqq \sum_{P\in \stLD(n)^{\bullet k}}(-1)^{\dinv(P)}t^{\area(P)}  \\
    S_{n,k} & \coloneqq \sum_{P\in \stLSQ(n)^{\bullet k}}(-1)^{\dinv(P)}t^{\area(P)}.
\end{align}
We will give a combinatorial formula for $S_{n,k}$ and explore the relationship between $S_{n,k}$ and $D_{n,k}$. We give here a table of these polynomials for small $n$ and $k$. 

\begin{minipage}{.9\textwidth}
    \centering
    \begin{tabular}{c||cc|cc|cc}
        $k$ & $S_{1,k}$ & $D_{1,k}$ & $S_{2,k}$ & $D_{2,k}$ & $S_{3,k}$ & $D_{3,k}$ \\
        \hline
        $0$&$1$&$1$&$0$&$t$&$t^{3} + t^{2} + t$&$t^{3} + t^{2}$\\ 
        $1$&&&$t + 1$&$1$&$0$&$t^{2} + 2 t$\\ 
        $2$&&&&&$t^{2} + t + 1$&$1$\\ 
    \end{tabular}
    \captionof{table}{Values of $S_{n,k}$ and $D_{n,k}$ for small $n,k$.} \label{tab:polvalues}      
\end{minipage}

\section{Schedule numbers}\label{sec:schedule}
In this paper, we will make extended use of \emph{schedule numbers}, a notion which allows us to factor the $q,t$-enumators of Dyck and square paths. They were introduced in \cite{Hicks2013} and have proven to be a very useful tool (see \cites{Sergel2017,HaglundSergel2021,IraciVandenWyngaerd2021}).
\begin{figure}
    \centering
    \begin{tikzpicture}[scale=.5]
        \draw[step=1.0, gray!60, thin] (0,0) grid (7,7);
        \draw[gray!60, thin] (1,0) -- (7,6);

        \draw[blue!60, line width=3pt] (0,0) -- (1,0) -- (1,1) -- (2,1) -- (2,2) -- (2,3) -- (3,3) -- (3,4) -- (3,5) -- (4,5) -- (5,5) -- (5,6) -- (6,6) -- (6,7) -- (7,7);

        \node at (1.5,0.5) {$1$};
        \draw (1.5,0.5) circle (.4cm);
        \node at (2.5,1.5) {$4$};
        \draw (2.5,1.5) circle (.4cm);
        \node at (2.5,2.5) {$5$};
        \draw (2.5,2.5) circle (.4cm);
        \node at (3.5,3.5) {$6$};
        \draw (3.5,3.5) circle (.4cm);
        \node at (3.5,4.5) {$7$};
        \draw (3.5,4.5) circle (.4cm);
        \node at (5.5,5.5) {$2$};
        \draw (5.5,5.5) circle (.4cm);
        \node at (6.5,6.5) {$3$};
        \draw (6.5,6.5) circle (.4cm);
        \node at (0--1-0.5,0+0.5) {$\bullet$};
        \node at (5-0-0.5,5+0.5) {$\bullet$};
        \node at (6-0-0.5,6+0.5) {$\bullet$};
    \end{tikzpicture}
    \caption{An element of $\stLSQ(7)^{\bullet 3}$.}\label{fig:stsquarepath}
\end{figure}
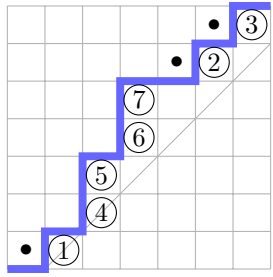

\begin{definition}
    Given $P\in \stLSQ(n)^{\bullet k}$ with shift $s$, set $\rho_i$ to be the labels of $P$ that lie in the $(i-s)$-th diagonal, written in decreasing order. Add a $\bullet$-symbol on top of each label that belongs to a decorated step. We define the \emph{diagonal word} of $P$ to be the word $\dw(P)\coloneqq \rho_0\rho_1\cdots$.
    The \emph{shifted diagonal word} of $P$ is the pair $\sdw(P)=(\dw(P), s)$.
\end{definition}
\begin{example}
    The diagonal word of the path in Figure~\ref{fig:stsquarepath} is $4\dec 165\dec 3\dec 2 7$.
\end{example}
The diagonal word of such a path is a \emph{decorated permutation.}
\begin{definition}
    Define the set of decorated permutations  $\mathfrak{S}_n^\bullet$ to be the set of permutations of $n$ where some letters are decorated with a $\bullet$. The subset where exactly $k$ letters are decorated is denoted by $\mathfrak{S}_n^{\bullet k}$.
\end{definition}
So for all $P\in \stLSQ(n)^{\bullet k}$, $\dw(P)\in\mathfrak{S}_n^{\bullet k}$. However, not every decorated permutation is the diagonal word of a path.

We recall the following classical definition.
\begin{definition}
    Given a permutation $\sigma$ of $n$, its \emph{major index} is \[\maj(\sigma) = \sum_{i: \sigma_i>\sigma_{i+1}}i.\] We denote by $\revmaj(\sigma)$ the major index of the reverse of $\sigma$ (i.e.\ the permutation $\sigma_n\sigma_{n-1}\cdots\sigma_{1}$).
\end{definition}

The following is easy to see.
\begin{proposition}
    All paths with a common diagonal word $\sigma$ share the same area, $\revmaj(\sigma)$.
\end{proposition}
Paths with the same shifted diagonal word must also share the same amount of bonus dinv.
\begin{definition}
    Given a decorated permutation $\sigma\in\mathfrak S_n^{\bullet}$ such that $\sigma = \rho_0 \rho_1\cdots\rho_l$ and the $\rho_i$ are the decreasing runs of $\sigma$. Given a shift $s\in\{0,\dots,l-1\}$, we call a run $\rho_i$ \emph{negative, zero} or \emph{positive} with respect to $s$, if $i<s$, $i=s$ or $i>s$, respectively.

    Let $\tilde \rho_i$ be the subword of undecorated letters of $\rho_i$. For $c \in \sigma$, define its \emph{(shifted) schedule number}\footnote{We will use the adjective ``shifted'' to emphasize that this value depends on the shift, but only when context requires it.} to be
    \[
        w_{\sigma,s}(c) =
        \begin{cases}
            \#\{d\in \tilde\rho_i \mid d>c\} + 1                                   & \text{if $\rho_i$ is zero and $c$ is undecorated}     \\
            \#\{d\in \tilde\rho_i \mid d>c\} + \#\{d\in \tilde\rho_{i-1}\mid d<c\} & \text{if $\rho_i$ is positive and $c$ is undecorated} \\
            \#\{d\in \tilde\rho_i \mid d<c\} + \#\{d\in \tilde\rho_{i+1}\mid d>c\} & \text{if $\rho_i$ is negative or $c$ is decorated} 
        \end{cases}
    \]
    For the sake of ease of notation, for any $s\in\N$ with $s\geq l$, set all the schedule numbers $w_{\sigma,s}(c)$ to be $0$.

    For a path $P\in\stLSQ(n)^{\bullet k}$ we define its \emph{schedule word} $\sched(P)$ to be the word whose letters are $w_{\sdw(P)}(c)$ for $c$ that runs through $\sdw(P)$.
\end{definition}
\begin{example} If we again take $\sigma$ to be the diagonal word of the path in Figure~\ref{fig:stsquarepath}, we compute
    \begin{align*}
        c\quad               & 4\dec 165\dec 3\dec 2 7 \\
        w_{\sigma,1}(c)\quad & 2212112
    \end{align*}
\end{example}

As in \cite{CorteelJosuatVergesVandenWyngaerd2023}, we will find the following definitions useful when working with schedule numbers.
\begin{definition}\label{def:cyclic-run}
    A \emph{cyclic (decreasing) run} of $\sigma\in \mathfrak{S_n}$ is a subword $\sigma_i\cdots \sigma_j$ such that there exists $m\in[n]$ such that 
    \[\sigma_i + m \pmod n,\dots, \sigma_j + m \pmod n\] is a decreasing sequence (where we take the representative modulo $n$ in $\{1,\dots, n\}$).
    
    Given $j\in[n]$ the \emph{left maximal cyclic run} of $\sigma_j$, $LMCR(\sigma_j)$ is defined as $\tau = \sigma_i\cdots\sigma_j$ where $i<j$ is minimal such that $\tau$ is a cyclic run. Similarly, given $i\in[n]$ the \emph{right maximal cyclic run} of $\sigma_i$, $RMCR(\sigma_j)$ is defined as $\tau = \sigma_i\cdots\sigma_j$ where $j>i$ is maximal such that $\tau$ is a cyclic run.
\end{definition}
See Figure~\ref{fig:cyclic-runs} for an illustration explaining the use of the terminology ``cyclic run''.
We can reformulate the definition of the shifted schedule word of a letter $\sigma_i$ to be
\[
    w_{\sigma,s}(\sigma_i) =
    \begin{cases}
        \#\{\sigma_k\in LMCR(\sigma_i) \mid k\neq i, \sigma_k \text{ undecorated}\} & \text{$\sigma_i$ undecorated and in a positive run} \\
        \#\{\sigma_k\in RMCR(\sigma_i)\mid k\neq i, \sigma_k \text{ undecorated}\}  & \text{$\sigma_i$ decorated or in negative run}
    \end{cases}
\]
Notice that the schedule numbers of undecorated letters in the zero run cannot be reformulated using cyclic runs.

The following is a specialization of \cite{IraciVandenWyngaerd2021}*{Theorem 5} to standardly labeled paths.
\begin{theorem}\label{thm:sched-formula} For all $n,s,k\in\N$ and $\sigma\in\mathfrak S_n^{\bullet k}$ we have
    \[\sum_{\substack{P\in\stLSQ(n)^{\bullet k}\\ \sdw(P) =(\sigma,s)}} q^{\dinv(P)}t^{\area(P)}= t^{\revmaj(\sigma)}q^{u(\sigma,s)}\prod_{c\in[n]}[w_s(c)]_q\] where $u(\sigma,s)$ is the number of undecorated letters in negative runs.
\end{theorem}
\begin{proof}[Sketch of proof]
    The general idea is to explicitly construct all the paths of a certain shifted diagonal  word $(\sigma,s)$ (i.e.\ with a fixed set of (decorated) labels in each diagonal). The strategy is to start with the empty path, and then insert, in all possible ways
    \begin{enumerate}
        \item  the undecorated labels into the $0$-diagonals,
        \item the undecorated labels into the $1$-diagonal, $2$-diagonal, etc.
        \item the undecorated labels into the $-1$-diagonal, then the $-2$-diagonal, etc.
    \end{enumerate}
    It turns out that the number of ways to insert a label $c$ into its diagonal is exactly the schedule number $w_{\sigma,s}(c)$, and that the contributions to the secondary and primary dinv is $q$-counted by $[w_{\sigma,s}(c)]_q$. The factor $u(\sigma,s)$ takes into account the bonus dinv.
    \begin{enumerate}
        \item[4.] Insert, in all possible ways, the decorated labels diagonal by diagonal (the order of the diagonals does not matter, since two decorated steps never form an attack relation).
    \end{enumerate}
    Here again, the  number of ways to insert is exactly $w_{\sigma,s}(c)$, and this time the entire dinv contribution is $q$-counted by $[w_{\sigma,s}(c)]_q$.
\end{proof}

From this theorem we readily deduce the following.

\begin{corollary}\label{cor:schedule}
    Given $(\sigma,s)\in \mathfrak S_n^{\bullet k}\times \N$
    \[\prod_{c\in[n]}w_{\sigma,s}(c)=\#\{P\in\stLSQ(n)^{\bullet k}\mid \sdw(P)=(\sigma,s)\}.
    \]
\end{corollary}

\section{Enumeration of decorated square paths at \texorpdfstring{$q=-1$}{q=-1}}\label{sec:enumeration}

When evaluating the schedule formula (Theorem~\ref{thm:sched-formula}) at $q=-1$, we notice that if $(\sigma,s)\in\mathfrak S_n^{\bullet}\times \N$ has an even schedule number, the sum over paths with shifted diagonal word $(\sigma,s)$ evaluates to $0$.  The following lemma tells us the same is true if $(\sigma,s)$ has at least one schedule $>1$.

\begin{lemma}\label{lem:scheds-are-interval}
    Let $(\sigma,s)\in\mathfrak S_n\times \N$ such that $w_{\sigma,s}(c)>0$ for all $c \in \sigma$. Then there exists a $j\in\N$ such that  $\{w_{\sigma,s}(c)\mid c\in \sigma\}=[j]$.
\end{lemma}
\begin{proof}
    The proof is similar to the one of the analogous statement for Dyck paths in \cite{CorteelJosuatVergesVandenWyngaerd2023}*{Lemma~3.14}. Denote by $\rho_0\rho_1\cdots$ the decreasing runs of $\sigma$. Let $\tau^{-}\coloneqq \rho_0\dots\rho_{s-1}, \tau^0\coloneqq \rho_s, \tau^{+} \coloneqq \rho_{s+1}\rho_{s+1}\cdots$ be the subwords of $\sigma$ consisting of the negative, zero and positive runs, respectively.
    Then, since the computation of schedules in zero and positive runs never depends on the negative runs, the $0$-shifted schedule numbers of the letters of $\tau^0\tau^{+}$ are exactly the $s$-shifted schedules of those letters in $\sigma$, and are thus $>0$. So we have \[\{w_{\tau^0\tau^{+},0}(c)\mid c\in \tau^0\tau^{+}\}=\{w_{\sigma,s}(c)\mid c\in \tau^0\tau^{+}\},\] and we can apply exactly the same argument as in \cite{CorteelJosuatVergesVandenWyngaerd2023}*{Lemma~3.14} to $\tau^0\tau^{+}$ to argue that this set is an interval of the form $[j_1]$ for some $j_1\in\N$.

    Now consider the elements of $\tau^{-}$. Here, the schedule of $c$ is computed by counting the number of undecorated elements different from $c$ in $c$'s right maximal cyclic run. Notice that such a right maximal cyclic run never extends beyond $\tau^0$. Set $l$ to be the number of letters in $\tau^0$ and take $c\in\tau^{-}$ such that $w_{\sigma,s}(c)>l$. We will show that there exists $d\in\tau^{-}$ such that $w_{\sigma,s}(d) = w_{\sigma,s}(c)-1$.
    Let $R=RMCR(c)$, take $c'$ be the leftmost letter in $R$, different from $c$ and $R' = RMCR(c')$. We must have that $R$ and $R'$ coincide on all undecorated letters of $R$, except $c$. So $w_{\sigma,s}(c')\geq w_{\sigma,s}(c)-1$. If the equality holds we are done, if not we iterate and look at $c''$ the rightmost letter in $R'$, different from $c'$. The leftmost letter in $\tau^{-}$ has $s$-shifted schedule number $\leq l$. Thus, since we took $w_{\sigma,s}(c)>l$, at some point in our iterations, we must find a letter whose $s$-shifted schedule is exactly $w_{\sigma,s}(c)-1$. It follows that \[\{w_{\sigma,s}(c)\mid c\in \tau^{-}\}=\{j_2,\dots,j_3\}\] for some $j_2,j_3\in\N$ with $j_2\leq l$. Since $\{w_{\sigma,s}(c)\mid c\in \tau^0\}=[l]$, the result follows.
\end{proof}

In this paper, we will give an enumerative formula for $S_{n,k}$. One can easily see that, combining Lemma~\ref{lem:scheds-are-interval} and Theorem~\ref{thm:sched-formula} we can conclude that when evaluating \eqref{eq:standard-delta-square} at $q=-1$, we are left only with paths whose schedule word is $1^n$:
\begin{equation}\label{eq:sched1sum}
    S_{n,k}=\sum_{\substack{P\in\stLSQ(n)^{\bullet k}\\ \sched(P)=1^n}}(-1)^{\dinv(P)}t^{\area(P)}.
\end{equation}
This formula is not yet ``enumerative'': the sum in the right-hand side of this equation still has some negative terms. In other words, there are still some cancellations to be accounted for.
In view of Corollary~\ref{cor:schedule}, we know that a path whose schedule word is $1^n$ is entirely determined by its shifted diagonal word $(\sigma,s)$. Furthermore, from Theorem~\ref{thm:sched-formula}, we know that the dinv of such a path is entirely accounted for by the number of undecorated letters in negative diagonals, denoted $u(\sigma,s)$.
Thus, we may reformulate \eqref{eq:sched1sum} as follows:
\begin{equation}
    S_{n,k}=\sum_{\substack{(\sigma,s) \in\mathfrak S_n^{\bullet k}\\ w_{\sigma,s}(\sigma_i)=1\;\forall i}}(-1)^{u(\sigma,s)}t^{\revmaj(\sigma)}.
\end{equation}
A given permutation $\sigma$ may contribute to the sum on the right-hand side multiple times: once for every shift for which its schedules are $1$. It turns out that all the cancellations happening in the sum on the right-hand side can be seen as to be leaving exactly one ``representative'' per distinct permutation $\sigma$. This motivates the following definition.
\begin{definition}\label{def:ADR}
    A permutation $\sigma \in \mathfrak{S}_n^{\bullet k}$ is an \emph{alternating dinv representative} (ADR) if, for some shift $s$, $(\sigma,s)$ has all (shifted) schedule numbers equal to $1$. It is called a \emph{Dyck} alternating dinv representative if this shift may be $0$. Denote the sets of such representatives by $\ADR_{n,k}$ and $\DADR_{n,k}$, respectively.
\end{definition}
\begin{example}
    For example, $\sigma = \dec 784\dec 23561$ is an \emph{alternating dinv representative}, since for shift $2$ and $3$ all the schedule numbers equal $1$. It is not a Dyck alternating dinv representative because $w_{\sigma,0}(8)=0$. See Example~\ref{ex:allshift-schedules} for more details.
\end{example}
In \cite{CorteelJosuatVergesVandenWyngaerd2023}, the authors established an enumeration formula for $D_{n,k}$ in terms of Dyck ADR's.
\begin{theorem} \label{thm:Dnk-enumeration}
    For all $n,k\in \N$, we have
    \[D_{n,k} = \sum_{\sigma \in \DADR_{n,k}}t^{\revmaj(\sigma)}\]
\end{theorem}
In this paper, we will show the following enumeration formula for $S_{n,k}$ in terms of ADR's.
\begin{theorem}[Cancellation Theorem -- word formulation]\label{thm:main1}
    For all $n,k\in \N$, we have
    \[S_{n,k}
        =\begin{cases}
            \sum\limits_{\sigma \in \ADR_{n,k}}t^{\revmaj(\sigma)} & \text{if $n-k$ is odd} \\
            0                                                      & \text{if $n-k$ is even.}
        \end{cases}\]
\end{theorem}
The proof of this theorem relies on a nice combinatorial explanation for the cancellations in \eqref{eq:sched1sum}, which is the subject of the next section.

\section{Cutting cycles}\label{sec:cutting}

\begin{figure}
    \centering
    \begin{subfigure}{.32\textwidth}
        \centering
        \begin{tikzpicture}[scale = .5]
            \draw[step=1.0, gray!60, thin] (0,0) grid (3,3);

            \draw[gray!60, thin] (0,0) -- (3,3);

            \draw[blue!60, line width=3pt] (0,0) -- (0,1) -- (0,2) -- (1,2) -- (2,2) -- (2,3) -- (3,3);

            \node at (0.5,0.5) {$1$};
            \draw (0.5,0.5) circle (.4cm);
            \node at (0.5,1.5) {$2$};
            \draw (0.5,1.5) circle (.4cm);
            \node at (2.5,2.5) {$3$};
            \draw (2.5,2.5) circle (.4cm);
            \node at (2-0-0.5,2+0.5) {$\bullet$};
        \end{tikzpicture}
        \subcaption{$P=\psi_3(P)$}
    \end{subfigure}
    \begin{subfigure}{.32\textwidth}
        \centering
        \begin{tikzpicture}[scale = .5]

            \draw[step=1.0, gray!60, thin] (0,0) grid (3,3);

            \draw[gray!60, thin] (1,0) -- (3,2);

            \draw[blue!60, line width=3pt] (0,0) -- (1,0) -- (1,1) -- (2,1) -- (2,2) -- (2,3) -- (3,3);

            \node at (1.5,0.5) {$3$};
            \draw (1.5,0.5) circle (.4cm);
            \node at (2.5,1.5) {$1$};
            \draw (2.5,1.5) circle (.4cm);
            \node at (2.5,2.5) {$2$};
            \draw (2.5,2.5) circle (.4cm);
            \node at (0--1-0.5,0+0.5) {$\bullet$};

        \end{tikzpicture}
        \subcaption{$\psi_1(P)$}
    \end{subfigure}
    \begin{subfigure}{.32\textwidth}
        \centering
        \begin{tikzpicture}[scale = .5]
            \draw[step=1.0, gray!60, thin] (0,0) grid (3,3);

            \draw[gray!60, thin] (0,0) -- (3,3);

            \draw[blue!60, line width=3pt] (0,0) -- (0,1) -- (1,1) -- (1,2) -- (1,3) -- (2,3) -- (3,3);

            \node at (0.5,0.5) {$3$};
            \draw (0.5,0.5) circle (.4cm);
            \node at (1.5,1.5) {$1$};
            \draw (1.5,1.5) circle (.4cm);
            \node at (1.5,2.5) {$2$};
            \draw (1.5,2.5) circle (.4cm);
            \node at (0-0-0.5,0+0.5) {$\bullet$};
        \end{tikzpicture}
        \caption{$\psi_2(P)$}
    \end{subfigure}
    \caption{An element of $\LSQ(3)^{\bullet 1}$ and its images by the $\psi_i$.}\label{fig:psi}
\end{figure}
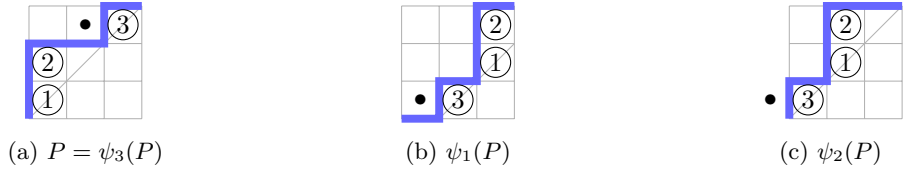

\begin{definition}
    Given $P\in\LSQ(n)^{\bullet k}$ and $i\in [n]$, split $P$ into two decorated labeled subpaths $P_1P_2$ where $P_1$ ends with the $i$-th east step of $P$. Define $\psi_i(P)$ to be $P_2P_1$. 
\end{definition}
We have $\psi_n(P) = P$. 
Furthermore, notice that $\psi_i(P)$ is not always a valid decorated labeled square path.
For example, in Figure~\ref{fig:psi}, we draw a path $P\in\LSQ(3)^{\bullet 1}$ and its images by the $\psi_i$. We see that $\psi_2(P)$ is not a valid decorated labeled square path because the first step is not a contractible valley, but is decorated. It is easy to see that if the step following the $i$-th horizontal step of $P\in\stLSQ(n)^{\bullet k}$ is not a decorated valley then $\psi_i(P)\in\LSQ(n)^{\bullet k}$.

\begin{definition}
    Given $P\in\LSQ(n)^{\bullet k}$ its \emph{cutting cycle} is \[\CC(P)\coloneqq \{\psi_i(P)\mid i\in [n]\}\cap \LSQ(n)^{\bullet k}.\]
    The name ``cutting cycle'' refers to the fact that each $Q\in \CC(P)$ is obtained by cutting $P$ into pieces and pasting those pieces back in a different order.
\end{definition}
In Figure~\ref{fig:psi}, the cutting cycle of $P$ thus consists of $\{P,\psi_1(P)\}$. Let us make the following easy observations.

\begin{observations}\label{obs:cc}\leavevmode
    \begin{enumerate}
        \item The cutting cycles partition $\LSQ(n)^{\bullet k}$.
        \item For all $P\in \LSQ(n)^{\bullet k}$, $\# \CC(P) = n-k$.
        \item Elements in the same cutting cycle have the same diagonal word.
              \item\label{it:area}Elements in the same cutting cycle have the same area.
    \end{enumerate}
\end{observations}

\begin{definition}
    We define an equivalence relation $\sim$ on the set $\LSQ(n)^{\bullet k}$ by setting $P\sim Q$ if and only if $Q\in \CC(P)$.
\end{definition}

\begin{definition}
    For $n,k\in\N$, we define the set of equivalence classes of paths with schedule word $1^n$:
    \[\mathcal S_{n,k} \coloneqq \{P\in \stLSQ(n)^{\bullet k}\mid \sched(P) = 1^n\}/\sim.\]
    For $C\in \mathcal S_{n,k}$ we denote by $\area(C)$ and $\dw(C)$ the common area and diagonal word of all the paths in the equivalence class $C$
\end{definition}
Keep in mind that not all the paths in the cutting cycle of a path of schedule word $1^n$ must have schedule word $1^n$ (see Example~\ref{ex:allshift-schedules}). Indeed, when the zero diagonal of a path contains more than one undecorated labels, its schedule is not $1^n$. This however, is the only restriction, as witnessed by the following result, which will follow from the closer study of the ADR object in the next section (Proposition~\ref{prop:decorating-algos}).

\begin{proposition}\label{prop:sched1-in-cc}
    Let $P\in\LSQ(n)^{\bullet k}$ such that $\sched(P) = 1^n$. Then for any $Q\in\mathcal{CC}(P)$, $\sched(Q) = 1^n$ if and only if $Q$'s main diagonal contains exactly one undecorated step. 
\end{proposition}

\begin{theorem}[Cancellation Theorem -- path formulation]\label{thm:main1.2}
    We have
    \[S_{n,k}
        =\begin{cases}
            \sum\limits_{C\in \mathcal S_{n,k}}t^{\area(C)} & \text{if $n-k$ is odd}  \\
            0                                               & \text{if $n-k$ is even}
        \end{cases}.\]
\end{theorem}

The remainder of this section will be dedicated to the proof this Theorem. We will also show that it is equivalent to Theorem~\ref{thm:main1}.

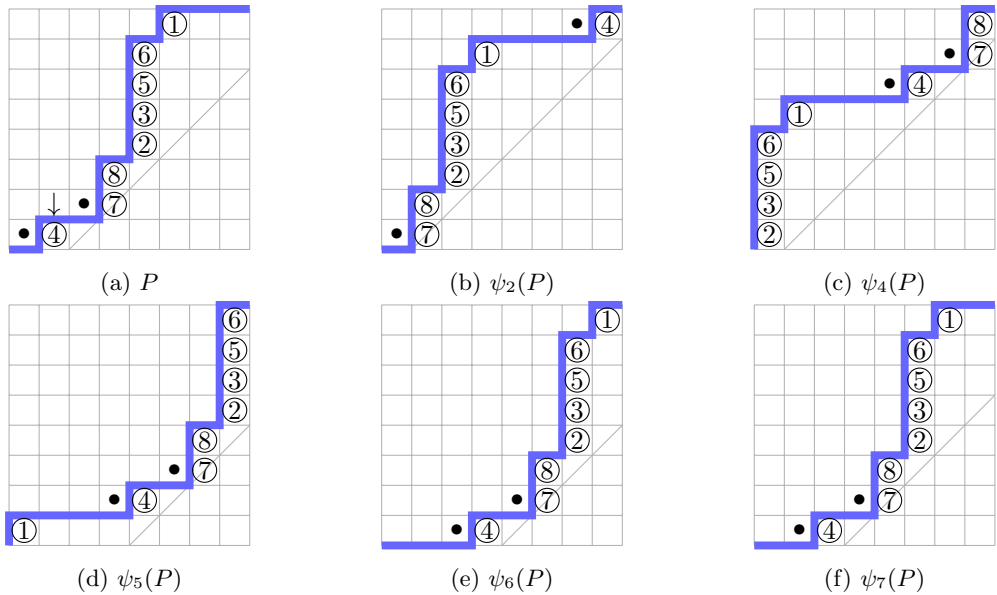
\begin{figure}
    \centering
    \begin{subfigure}{.32\textwidth}
        \centering
        \begin{tikzpicture}[scale=.4]
            \draw[step=1.0, gray!60, thin] (0,0) grid (8,8);
            \draw[gray!60, thin] (2,0) -- (8,6);

            \draw[blue!60, line width=3pt] (0,0) -- (1,0) -- (1,1) -- (2,1) -- (3,1) -- (3,2) -- (3,3) -- (4,3) -- (4,4) -- (4,5) -- (4,6) -- (4,7) -- (5,7) -- (5,8) -- (6,8) -- (7,8) -- (8,8);

            \node at (1.5,0.5) {$4$};
            \draw (1.5,0.5) circle (.4cm);
            \node at (3.5,1.5) {$7$};
            \draw (3.5,1.5) circle (.4cm);
            \node at (3.5,2.5) {$8$};
            \draw (3.5,2.5) circle (.4cm);
            \node at (4.5,3.5) {$2$};
            \draw (4.5,3.5) circle (.4cm);
            \node at (4.5,4.5) {$3$};
            \draw (4.5,4.5) circle (.4cm);
            \node at (4.5,5.5) {$5$};
            \draw (4.5,5.5) circle (.4cm);
            \node at (4.5,6.5) {$6$};
            \draw (4.5,6.5) circle (.4cm);
            \node at (5.5,7.5) {$1$};
            \draw (5.5,7.5) circle (.4cm);
            \node at (0--1-0.5,0+0.5) {$\bullet$};
            \node at (1--2-0.5,1+0.5) {$\bullet$};

            \node at (1.5,1.5)  {$\downarrow$};
        \end{tikzpicture}
        \subcaption{$P$}
        \label{fig:cc-P}
    \end{subfigure}
    \begin{subfigure}{.32\textwidth}
        \centering
        \begin{tikzpicture}[scale=.4]
            \draw[step=1.0, gray!60, thin] (0,0) grid (8,8);

            \draw[gray!60, thin] (1,0) -- (8,7);

            \draw[blue!60, line width=3pt] (0,0) -- (1,0) -- (1,1) -- (1,2) -- (2,2) -- (2,3) -- (2,4) -- (2,5) -- (2,6) -- (3,6) -- (3,7) -- (4,7) -- (5,7) -- (6,7) -- (7,7) -- (7,8) -- (8,8);

            \node at (1.5,0.5) {$7$};
            \draw (1.5,0.5) circle (.4cm);
            \node at (1.5,1.5) {$8$};
            \draw (1.5,1.5) circle (.4cm);
            \node at (2.5,2.5) {$2$};
            \draw (2.5,2.5) circle (.4cm);
            \node at (2.5,3.5) {$3$};
            \draw (2.5,3.5) circle (.4cm);
            \node at (2.5,4.5) {$5$};
            \draw (2.5,4.5) circle (.4cm);
            \node at (2.5,5.5) {$6$};
            \draw (2.5,5.5) circle (.4cm);
            \node at (3.5,6.5) {$1$};
            \draw (3.5,6.5) circle (.4cm);
            \node at (7.5,7.5) {$4$};
            \draw (7.5,7.5) circle (.4cm);
            \node at (0--1-0.5,0+0.5) {$\bullet$};
            \node at (7-0-0.5,7+0.5) {$\bullet$};
        \end{tikzpicture}
        \subcaption{$\psi_2(P)$}
        \label{fig:cc-d0}
    \end{subfigure}
    \begin{subfigure}{.32\textwidth}
        \centering
        \begin{tikzpicture}[scale=.4]
            \draw[step=1.0, gray!60, thin] (0,0) grid (8,8);

            \draw[gray!60, thin] (1,0) -- (8,7);

            \draw[blue!60, line width=3pt] (0,0) -- (0,1) -- (0,2) -- (0,3) -- (0,4) -- (1,4) -- (1,5) -- (2,5) -- (3,5) -- (4,5) -- (5,5) -- (5,6) -- (6,6) -- (7,6) -- (7,7) -- (7,8) -- (8,8);

            \node at (0.5,0.5) {$2$};
            \draw (0.5,0.5) circle (.4cm);
            \node at (0.5,1.5) {$3$};
            \draw (0.5,1.5) circle (.4cm);
            \node at (0.5,2.5) {$5$};
            \draw (0.5,2.5) circle (.4cm);
            \node at (0.5,3.5) {$6$};
            \draw (0.5,3.5) circle (.4cm);
            \node at (1.5,4.5) {$1$};
            \draw (1.5,4.5) circle (.4cm);
            \node at (5.5,5.5) {$4$};
            \draw (5.5,5.5) circle (.4cm);
            \node at (7.5,6.5) {$7$};
            \draw (7.5,6.5) circle (.4cm);
            \node at (7.5,7.5) {$8$};
            \draw (7.5,7.5) circle (.4cm);
            \node at (5-0-0.5,5+0.5) {$\bullet$};
            \node at (6--1-0.5,6+0.5) {$\bullet$};

        \end{tikzpicture}
        \subcaption{$\psi_4(P)$}
        \label{fig:cc-d1}
    \end{subfigure}
    \begin{subfigure}{.32\textwidth}
        \centering
        \begin{tikzpicture}[scale=.4]
            \draw[step=1.0, gray!60, thin] (0,0) grid (8,8);

            \draw[gray!60, thin] (4,0) -- (8,4);

            \draw[blue!60, line width=3pt] (0,0) -- (0,1) -- (1,1) -- (2,1) -- (3,1) -- (4,1) -- (4,2) -- (5,2) -- (6,2) -- (6,3) -- (6,4) -- (7,4) -- (7,5) -- (7,6) -- (7,7) -- (7,8) -- (8,8);

            \node at (0.5,0.5) {$1$};
            \draw (0.5,0.5) circle (.4cm);
            \node at (4.5,1.5) {$4$};
            \draw (4.5,1.5) circle (.4cm);
            \node at (6.5,2.5) {$7$};
            \draw (6.5,2.5) circle (.4cm);
            \node at (6.5,3.5) {$8$};
            \draw (6.5,3.5) circle (.4cm);
            \node at (7.5,4.5) {$2$};
            \draw (7.5,4.5) circle (.4cm);
            \node at (7.5,5.5) {$3$};
            \draw (7.5,5.5) circle (.4cm);
            \node at (7.5,6.5) {$5$};
            \draw (7.5,6.5) circle (.4cm);
            \node at (7.5,7.5) {$6$};
            \draw (7.5,7.5) circle (.4cm);
            \node at (1--3-0.5,1+0.5) {$\bullet$};
            \node at (2--4-0.5,2+0.5) {$\bullet$};

        \end{tikzpicture}
        \subcaption{$\psi_5(P)$}
        \label{fig:cc-d5}
    \end{subfigure}
    \begin{subfigure}{.32\textwidth}
        \centering
        \begin{tikzpicture}[scale=.4]
            \draw[step=1.0, gray!60, thin] (0,0) grid (8,8);

            \draw[gray!60, thin] (4,0) -- (8,4);

            \draw[blue!60, line width=3pt] (0,0) -- (1,0) -- (2,0) -- (3,0) -- (3,1) -- (4,1) -- (5,1) -- (5,2) -- (5,3) -- (6,3) -- (6,4) -- (6,5) -- (6,6) -- (6,7) -- (7,7) -- (7,8) -- (8,8);

            \node at (3.5,0.5) {$4$};
            \draw (3.5,0.5) circle (.4cm);
            \node at (5.5,1.5) {$7$};
            \draw (5.5,1.5) circle (.4cm);
            \node at (5.5,2.5) {$8$};
            \draw (5.5,2.5) circle (.4cm);
            \node at (6.5,3.5) {$2$};
            \draw (6.5,3.5) circle (.4cm);
            \node at (6.5,4.5) {$3$};
            \draw (6.5,4.5) circle (.4cm);
            \node at (6.5,5.5) {$5$};
            \draw (6.5,5.5) circle (.4cm);
            \node at (6.5,6.5) {$6$};
            \draw (6.5,6.5) circle (.4cm);
            \node at (7.5,7.5) {$1$};
            \draw (7.5,7.5) circle (.4cm);
            \node at (0--3-0.5,0+0.5) {$\bullet$};
            \node at (1--4-0.5,1+0.5) {$\bullet$};

        \end{tikzpicture}
        \subcaption{$\psi_6(P)$}
        \label{fig:cc-d4}
    \end{subfigure}
    \begin{subfigure}{.32\textwidth}
        \centering
        \begin{tikzpicture}[scale=.4]
            \draw[step=1.0, gray!60, thin] (0,0) grid (8,8);

            \draw[gray!60, thin] (3,0) -- (8,5);

            \draw[blue!60, line width=3pt] (0,0) -- (1,0) -- (2,0) -- (2,1) -- (3,1) -- (4,1) -- (4,2) -- (4,3) -- (5,3) -- (5,4) -- (5,5) -- (5,6) -- (5,7) -- (6,7) -- (6,8) -- (7,8) -- (8,8);

            \node at (2.5,0.5) {$4$};
            \draw (2.5,0.5) circle (.4cm);
            \node at (4.5,1.5) {$7$};
            \draw (4.5,1.5) circle (.4cm);
            \node at (4.5,2.5) {$8$};
            \draw (4.5,2.5) circle (.4cm);
            \node at (5.5,3.5) {$2$};
            \draw (5.5,3.5) circle (.4cm);
            \node at (5.5,4.5) {$3$};
            \draw (5.5,4.5) circle (.4cm);
            \node at (5.5,5.5) {$5$};
            \draw (5.5,5.5) circle (.4cm);
            \node at (5.5,6.5) {$6$};
            \draw (5.5,6.5) circle (.4cm);
            \node at (6.5,7.5) {$1$};
            \draw (6.5,7.5) circle (.4cm);
            \node at (0--2-0.5,0+0.5) {$\bullet$};
            \node at (1--3-0.5,1+0.5) {$\bullet$};
        \end{tikzpicture}
        \subcaption{$\psi_7(P)$}
        \label{fig:cc-d3}
    \end{subfigure}
    \caption{The equivalence class of a path $P$ of schedule word $1^n$}
    \label{fig:cutting-cycle-1n}
\end{figure}

\begin{example}\label{ex:allshift-schedules}
    In Figure~\ref{fig:cutting-cycle-1n}, we have displayed all the elements of the cutting cycle of a schedule $1^n$ path $P$, i.e.\ an element $\mathcal S_{n,k}$, for $n=8$ and $k=2$. Here we compute the schedule numbers for all possible shifts: there are $5$ runs in the common diagonal word $\dec 784\dec 23561$, so the nonzero schedules occur for a shift between $0$ and $4$:
    \begin{align*}
        c\quad               & \dec 784\dec 23561 \\
        w_{\sigma,0}(c)\quad & 10111111           \\
        w_{\sigma,1}(c)\quad & 11121111           \\
        w_{\sigma,2}(c)\quad & 11111111           \\
        w_{\sigma,3}(c)\quad & 11111111           \\
        w_{\sigma,4}(c)\quad & 11111112.
    \end{align*}
    We confirm that $P$ (see Figure~\ref{fig:cc-P}) has shift $2$ and schedule word $1^n$. There are no paths with diagonal word $\dec 784\dec 23561$ and shift $0$, since one of the schedule numbers is 0. In $P$'s cutting cycle there are $2$ paths of shift $1$, $1$ of shift $3$ and $2$ of shift $4$.
\end{example}

For each equivalence class in $\mathcal S_{n,k}$, we fix a canonical representative.
\begin{definition}\label{def:Ptilde} Given $P \in \LSQ(n)^{\bullet k}$, its \emph{breaking point} is:
    \begin{enumerate}
        \item The first (i.e., leftmost) starting point of an undecorated vertical step on the bottom diagonal;
        \item if no such steps exist, the first starting point of a decorated vertical step on the bottom diagonal.
    \end{enumerate}

    In the first case, the \emph{breaking step} is the horizontal step immediately to the left of the breaking point. In the second case the breaking point is either preceded by two consecutive east steps, or the breaking point is the endpoint of the first horizontal step of $P$ and the last step of $P$ is a horizontal step by definition. Thus, in the second case, the \emph{breaking step} is the horizontal step two to the left (read cyclically) of the breaking point.
    
    If the breaking step is the $i$th horizontal step, we define $\tilde{P} \coloneqq \psi_i(P)$.
\end{definition}
\begin{example}\label{ex:Ptilde}
    The breaking step of $P$ in Figure~\ref{fig:cc-P} is marked with an arrow. It is the second horizontal step of $P$, thus, the canonical representative of $P$'s cutting cycle is $\psi_2(P)$, depicted in Figure~\ref{fig:cc-d0}.
\end{example}

\begin{figure}
    \centering
    \begin{tikzpicture}[scale = .35,every node/.style={scale=.8}]
        \draw[step=1.0, gray!60, thin] (0,0) grid (24,24);

        \draw[gray!60, thin] (3,0) -- (24,21);

        \draw[ultra thick, dashed] (0,0) rectangle (7,4) rectangle (11,20) rectangle (24,24);
        \draw[blue!60, line width=2pt] (0,0) -- (1,0) -- (1,1) -- (2,1) -- (3,1) -- (3,2) -- (4,2) -- (4,3) -- (5,3) -- (6,3) -- (6,4) -- (7,4) -- (7,5) -- (7,6) -- (7,7) -- (8,7) -- (8,8) -- (8,9) -- (8,10) -- (8,11) -- (8,12) -- (8,13) -- (8,14) -- (9,14) -- (9,15) -- (9,16) -- (10,16) -- (10,17) -- (10,18) -- (10,19) -- (10,20) -- (11,20) -- (12,20) -- (13,20) -- (14,20) -- (14,21) -- (15,21) -- (16,21) -- (17,21) -- (18,21) -- (19,21) -- (19,22) -- (20,22) -- (20,23) -- (21,23) -- (22,23) -- (23,23) -- (23,24) -- (24,24);

        \node at (1.5,0.5) {$9$};
        \draw (1.5,0.5) circle (.4cm);
        \node at (3.5,1.5) {$12$};
        \draw (3.5,1.5) circle (.4cm);
        \node at (4.5,2.5) {$14$};
        \draw (4.5,2.5) circle (.4cm);
        \node at (6.5,3.5) {$2$};
        \draw (6.5,3.5) circle (.4cm);
        \node at (7.5,4.5) {$3$};
        \draw (7.5,4.5) circle (.4cm);
        \node at (7.5,5.5) {$7$};
        \draw (7.5,5.5) circle (.4cm);
        \node at (7.5,6.5) {$10$};
        \draw (7.5,6.5) circle (.4cm);
        \node at (8.5,7.5) {$4$};
        \draw (8.5,7.5) circle (.4cm);
        \node at (8.5,8.5) {$6$};
        \draw (8.5,8.5) circle (.4cm);
        \node at (8.5,9.5) {$8$};
        \draw (8.5,9.5) circle (.4cm);
        \node at (8.5,10.5) {$11$};
        \draw (8.5,10.5) circle (.4cm);
        \node at (8.5,11.5) {$13$};
        \draw (8.5,11.5) circle (.4cm);
        \node at (8.5,12.5) {$19$};
        \draw (8.5,12.5) circle (.4cm);
        \node at (8.5,13.5) {$22$};
        \draw (8.5,13.5) circle (.4cm);
        \node at (9.5,14.5) {$17$};
        \draw (9.5,14.5) circle (.4cm);
        \node at (9.5,15.5) {$21$};
        \draw (9.5,15.5) circle (.4cm);
        \node at (10.5,16.5) {$16$};
        \draw (10.5,16.5) circle (.4cm);
        \node at (10.5,17.5) {$18$};
        \draw (10.5,17.5) circle (.4cm);
        \node at (10.5,18.5) {$23$};
        \draw (10.5,18.5) circle (.4cm);
        \node at (10.5,19.5) {$24$};
        \draw (10.5,19.5) circle (.4cm);
        \node at (14.5,20.5) {$1$};
        \draw (14.5,20.5) circle (.4cm);
        \node at (19.5,21.5) {$15$};
        \draw (19.5,21.5) circle (.4cm);
        \node at (20.5,22.5) {$20$};
        \draw (20.5,22.5) circle (.4cm);
        \node at (23.5,23.5) {$5$};
        \draw (23.5,23.5) circle (.4cm);
        \node at (0--1-0.5,0+0.5) {$\bullet$};
        \node at (1--2-0.5,1+0.5) {$\bullet$};
        \node at (2--2-0.5,2+0.5) {$\bullet$};
        \node at (3--3-0.5,3+0.5) {$\bullet$};
        \node at (20-6-0.5,20+0.5) {$\bullet$};
        \node at (21-2-0.5,21+0.5) {$\bullet$};
        \node at (22-2-0.5,22+0.5) {$\bullet$};
        \node at (23-0-0.5,23+0.5) {$\bullet$};
    \end{tikzpicture}
    \caption{A path $P$ with $\sched(P)=1^n$, with its three stretches as described in Lemma~\ref{lem:shape}}\label{fig:sched1path}
\end{figure}
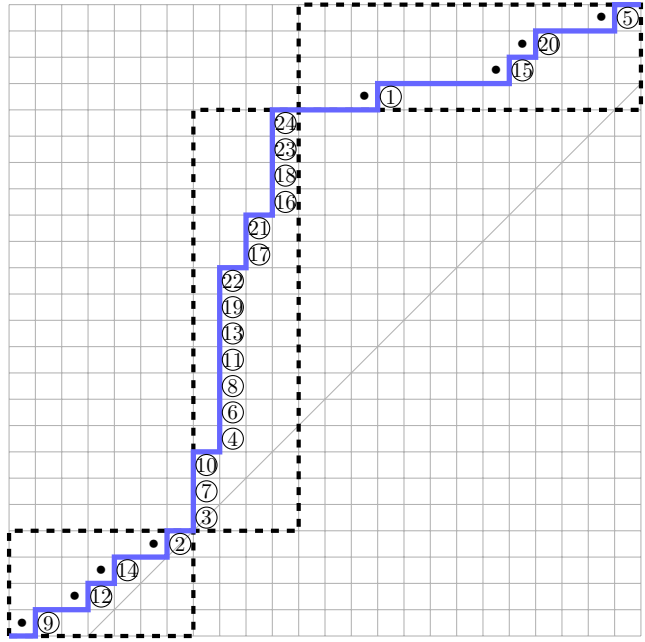

The condition for a path to have schedules all $1$ is quite restrictive. It will be useful to have a general idea of the shape of such paths. See Figure~\ref{fig:sched1path} for a large example of such a path.

\begin{lemma}[Shape lemma]\label{lem:shape}
    Let $P \in \LSQ(n)^{\bullet k}$ satisfy $\sched(P)=1^n$. Then $P$ consists of three stretches (indicated by the dashed rectangles in Figure~\ref{fig:sched1path}):
        \begin{enumerate}
            \item In the first stretch, all vertical steps are decorated contractible valleys (though there might be none), all of which lie in a negative diagonal.
            \item The second stretch, starts with the first undecorated vertical step of the path and ends with the horizontal step following the last undecorated vertical step of the path. All the vertical steps in this stretch are undecorated, and there are no two consecutive horizontal steps. Two labels in this stretch never form an attack relation.
            \item In the final stretch, all vertical steps are decorated contractible valleys (though there might be none). All labels lie in positive diagonals.
        \end{enumerate}
    \end{lemma}

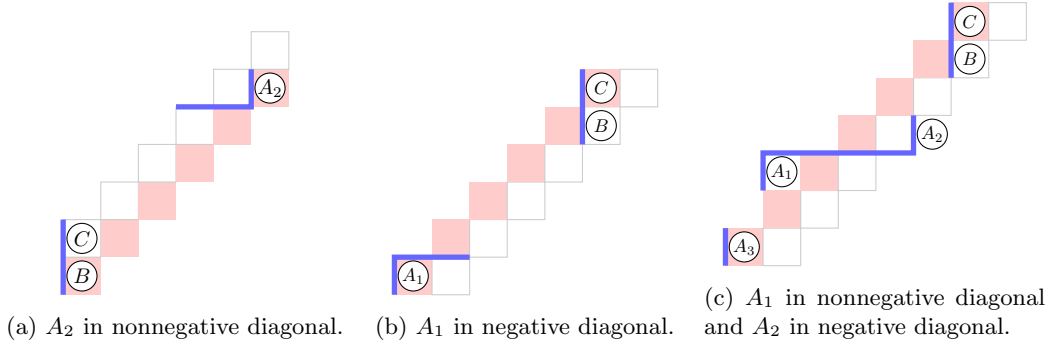
\begin{figure}
    \centering
    \begin{subfigure}{.3\textwidth}
        \centering
        \begin{tikzpicture}[scale = .5,every node/.style={scale=.8}]
            \foreach \i in {0,...,5}{
                \fill[red!20] (\i,\i) rectangle (\i+1,\i+1);
                    \draw[black!20] (\i,\i+1) rectangle (\i+1,\i+2);
                }
            \draw[blue!60, line width=2pt] (0,0)-- (0,2) (3,5)-|(5,6);
            \draw[fill = white] (5.5,5.5) circle(.4cm) node {$A_2$};
            \draw[fill = white] (.5,.5) circle(.4cm) node {$B$};
            \draw[fill = white] (.5,1.5) circle(.4cm) node {$C$};
        \end{tikzpicture}
        \caption{$A_2$ in nonnegative diagonal.}\label{fig:A2nonneg}
    \end{subfigure}
    \begin{subfigure}{.3\textwidth}
        \centering
        \begin{tikzpicture}[scale = .5,every node/.style={scale=.7}]
            \foreach \i in {0,...,5}{
                \fill[red!20] (\i,\i) rectangle (\i+1,\i+1);
                \draw[black!20] (\i+1,\i) rectangle (\i+2,\i+1);
                }
            \draw[blue!60, line width=2pt] (0,0)|-(2,1) (5,4)--(5,6);
            \draw[fill = white] (.5,.5) circle(.4cm) node {$A_1$};
            \draw[fill = white] (5.5,4.5) circle(.4cm) node {$B$};
            \draw[fill = white] (5.5,5.5) circle(.4cm) node {$C$};
        \end{tikzpicture}
        \caption{$A_1$ in negative diagonal.}\label{fig:A1neg}
    \end{subfigure}
    \begin{subfigure}{.3\textwidth}
        \centering
        \begin{tikzpicture}[scale = .5,every node/.style={scale=.7}]
            \foreach \i in {0,...,6}{
                \fill[red!20] (\i,\i) rectangle (\i+1,\i+1);
                \draw[black!20] (\i+1,\i) rectangle (\i+2,\i+1);
                }
                \coordinate (1) at (1,2);
                \coordinate (2) at (5,3);
                \draw[blue!60, line width=2pt] (1) |- (2) -- ++ (0,1) (6,5)--(6,7) (0,0)--(0,1);
            \draw[fill = white] (1) ++ (.5,.5) circle (.4cm) node {$A_1$};
            \draw[fill = white] (2) ++ (.5,.5) circle (.4cm) node {$A_2$};
            \draw[fill = white] (6.5,5.5) circle (.4cm) node {$B$};
            \draw[fill = white] (6.5,6.5) circle (.4cm) node {$C$};
            \draw[fill = white] (.5,.5) circle (.4cm) node {$A_3$};
        \end{tikzpicture}
        \caption{$A_1$ in nonnegative diagonal and $A_2$ in negative diagonal.}\label{fig:A1nonnegA2neg}
    \end{subfigure}
    \caption{Forced attack relations.}\label{fig:forced-attack}
\end{figure}

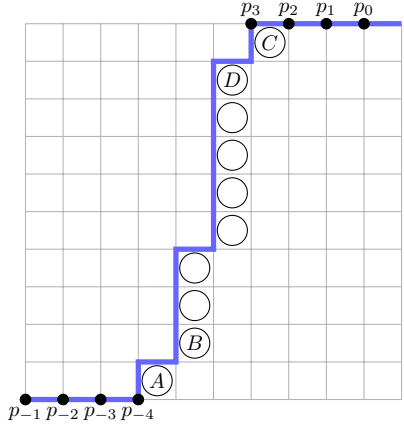
\begin{figure}
    \centering
    \begin{tikzpicture}[scale = .5,every node/.style={scale=.8}]
        \draw[gray!60,thin] (0,0) grid (10,10);
        \draw[blue!60, line width=2pt](0,0)-|(3,1)-|(4,4)-|(5,9)-|(6,10) --(10,10);
        \filldraw (9,10) circle(4pt) node[above] {$p_0$};
        \filldraw (8,10) circle(4pt) node[above] {$p_1$};
        \filldraw (7,10) circle(4pt) node[above] {$p_2$};
        \filldraw (6,10) circle(4pt) node[above] {$p_3$};
        \filldraw (0,0) circle(4pt) node[below] {$p_{-1}$};
        \filldraw (1,0) circle(4pt) node[below] {$p_{-2}$};
        \filldraw (2,0) circle(4pt) node[below] {$p_{-3}$};
        \filldraw (3,0) circle(4pt) node[below] {$p_{-4}$};
        \draw (3.5,.5) circle(.4cm) node {$A$};
        \draw (4.5,1.5) circle(.4cm) node {$B$};
        \draw (6.5,9.5) circle(.4cm) node {$C$};
        \draw (5.5,8.5) circle(.4cm) node {$D$};
        \draw (5.5,7.5) circle(.4cm) node {};
        \draw (5.5,6.5) circle(.4cm) node {};
        \draw (5.5,5.5) circle(.4cm) node {};
        \draw (5.5,4.5) circle(.4cm) node {};
        \draw (4.5,3.5) circle(.4cm) node {};
        \draw (4.5,2.5) circle(.4cm) node {};
    \end{tikzpicture}
    \caption{Insertion points for decorated valleys.}\label{fig:insertion-points}
\end{figure}

\begin{proof}
    The proof of this lemma relies heavily on the ideas in the proof of the schedule formula (Theorem~\ref{thm:sched-formula}), of which we gave a sketch above. Indeed, consider the construction algorithm of the unique path $P$ of a certain shifted diagonal word with schedules all $1$. If we stop at step 3, before inserting any of the labels attached to a decorated step, we must have obtained an undecorated path \emph{with no attack relations}, call it $P_1$ and denote its size by $n_1$.
    
    We claim that in $P_1$ all consecutive pairs of horizontal steps occur either before all or after all vertical steps of the path (in other words $P_1$ is of the form $E^*(NN^*E)^*E^*$). Suppose to the contrary that there exists a stretch of consecutive horizontal steps of size at least two that does not lie on the line $y=0$ or $y=n_1$. Let $A_1$ (respectively $A_2$) be the label of the vertical step immediately preceding (respectively succeeding) this horizontal stretch. We distinguish three cases.
    \begin{enumerate}
        \item $A_2$ lies in a nonnegative diagonal (see Figure~\ref{fig:A2nonneg}). We have that $A_2$'s vertical step is preceded by two horizontal steps that lie in diagonals say $i$ and $i+1$ with $i\geq 0$. Since the path starts at diagonal $0$, there must be, before the step labeled $A_2$, two consecutive vertical steps in diagonal $i$ and $i+1$, let them be labeled $B$ and $C$, respectively. If $C>A_2$ then their steps are in an attack relation. If $C<A_2$, then $B<A_2$ since $B<C$, and so the steps of $B$ and $A_2$ are in an attack relation.
        \item $A_1$ lies in a negative diagonal (see Figure~\ref{fig:A1neg}). Now $A_1$'s step is followed by two horizontal steps, in diagonals $i$ and $i-1$, with $i<0$. Since the path must end up strictly above the diagonal $x=y$, there must be consecutive vertical steps after $A_1$'s step, in diagonals $i-1$ and $i$, labeled $B$ and $C$, respectively. By similar arguments as before, there must be at least one attack relation between the steps labeled $A_1$, $B$ and $C$.
        \item $A_1$ lies in a nonnegative diagonal and $A_2$ lies in a negative diagonal (see Figure~\ref{fig:A1nonnegA2neg}). This means that there must be two consecutive horizontal steps in diagonals $0$ and $-1$. Since the path must end up strictly above $x=y$, there must be two consecutive vertical steps weakly after $A_2$'s step, in diagonals $-1$ and $0$, labeled $B$ and $C$ respectively (note that we might have $A_2=B$). Since $A_1$ is in a nonnegative diagonal and the path starts at $(0,0)$, there must be an $A_3$ in the $0$-diagonal, where we might have $A_3=A_1$. Between the steps labeled $A_3$, $B$ and $C$, there is at least one attack relation.
    \end{enumerate}
    Since $P_1$ has no attack relations, we have established the claim. Furthermore, in $P_1$ no diagonal contains more that $2$ (non-decorated) labels. Indeed, it is easy to see from the schedule number definition that if all schedules are $1$, each run of the diagonal word must contain at most two non-decorated letters.

    Let us now continue the construction algorithm in the proof of the schedule formula (Theorem~\ref{thm:sched-formula}). In step 4, the decorated labels are inserted into their diagonals. Since the schedule numbers are all $1$, for each decorated label, there must be exactly one valid way to do this. Let us observe that each run of the diagonal word of a path, except the first one corresponding to the bottom most diagonal of the path, must contain at least one non-decorated label (indeed all except the bottom most diagonal must contain a vertical step preceded by another vertical step, which is not a contractible valley). Thus, decorated valleys must be inserted into the already populated diagonals of $P_1$, plus maybe one diagonal below its bottom most diagonal. We will identify \emph{insertion points} on $P_1$, one for each diagonal into which we might insert decorated valleys. We will argue that inserting a horizontal step followed by a decorated vertical step with the appropriate label at the insertion points will yield a valid path. Keeping the insertion point to the left of the newly inserted steps, we may insert more than one decorated label, from biggest to smallest label, so that the consecutive decorated valleys are contractible.
    Since the way of inserting must be unique, we will obtain the unique path $P$ after these insertions.

    We denote by $p_i$ the insertion point of the $i$-th diagonal (we refer to Figure~\ref{fig:insertion-points} for a running example). Denote by $M$ (respectively $m$) the number of positive (respectively, negative) diagonals intersecting $P_1$ and recall that $n_1$ is the size of $P_1$. For $i=0,1,2,\dots, M$ the $p_i$ lie on the last horizontal stretch of $P_1$, from right to left: $p_0=(n_1-1,n_1), p_1=(n_1-2,n_1-1), \dots, p_M = (n_1-M-1,n_1)$.
    For $i=-1,-2,\dots,p_{-m-1}$ the $p_i$ lie on the first horizontal stretch of $P_1$ form left to right: $p_{-1}=(0,0), p_{-2} = (1,0), \dots, p_{-m-1} = (-m,0)$. It is immediately clear that for $i\neq M$ and $i\neq -m-1$, the insertion of decorated valleys as described yields a valid path (i.e.\ the decorated valleys are contractible and labels are increasing in columns). 
    
    When we insert a decorated valley $V$ at $p_{-m-1}$, $V$ ends up under a label in the $(-m)$-diagonal, say $A$, so we must check that $V<A$. If $A$ is the only undecorated label in its diagonal, then if $A<V$, the schedule of $V$ is $0$. If there are two undecorated labels in the $(-m)$-diagonal, say $A$ and $B$ from bottom to top, then since $P_1$ is attack-free, $A>B$, so if $A<V$, then $V$ is larger than both $A$ and $B$, and $V$'s schedule number equals 2. 
    
    Similarly, when inserting a decorated valley $V$ at $p_M$, $V$'s vertical step will be preceded by a single horizontal step that is preceded by the vertical step labeled say $C$. We must check that $V$'s valley is contractible, i.e.\ that $C<V$. If $C$ is the only decorated step in its diagonal, then if $C>V$, $V$'s schedule is zero. If there are two undecorated steps in the $M$-diagonal, say $C$ and $D$ bottom to top, we have $D>C$, so if $C>V$ then $V$ is smaller that both $C$ and $C$ and its schedule number is 2.

    It follows that the unique path $P$ is obtained from $P_1$ by inserting decorated valleys into its first and last horizontal stretch as described. The general shape of the path described in the statement of the lemma is a consequence of this fact, and the previously discussed shape of $P_1$.
\end{proof}

Using this result and recalling Definition~\ref{def:Ptilde}, we can deduce the shape of canonical representatives of schedule $1^n$ paths.
\begin{corollary}\label{cor:shape-canonical}
    Let $P\in\stLSQ(n)^{\bullet k}$ such that $\sched(P) = 1^n$. Then $\tilde P$ consists of three stretches satisfying the same shape restrictions as described in Lemma~\ref{lem:shape}. Furthermore, all the labels in the first stretch of the path must lie in the $(-1)$-diagonal. 
\end{corollary}

\begin{lemma}\label{PTildedinv0}
    For all $P \in \stLSQ(n)$ such that $\sched(P)=1^n$, $\dinv(\widetilde{P}) = 0$.
\end{lemma}
\begin{proof}
    We write $P = P_1P_2$, where $P_1$'s last step is the breaking step. Thus, we have $\widetilde P = P_2P_1$.
    Let us show that $\dinv(\tilde P)=0$.
    Since $\sched(P)=1^n$ it has no attack relations (Observations~\ref{obs:cc}). It follows that there are no attack relations among two labels in $P_1$ or two labels in $P_2$.
    Furthermore, by the choice of the breaking step, the only labels in a negative diagonal are attached to decorated steps in the $(-1)$-diagonal. We can think of the tertiary dinv of these labels as being ``cancelled out'' by the $-1$ each of these decorations contributes to the dinv.
    So all we need to consider is the attack relations between a label in $P_1$ and one in $P_2$. By Lemma~\ref{lem:shape} we know that all the up steps in $P_1$ are decorated. Let $a$ be a label of $P_1$. Since the schedule number of $a$ is $1$, there is exactly one undecorated label $b$ in $P$ that is smaller and in $a$'s diagonal or bigger and in the diagonal above $a$'s. These are the only possible attach relations in $\tilde P$. Since $b$ is undecorated, $b$ must be a label of $P_2$. So in $P$, $a$'s step and $b$'s step are not in an attack relation, but in $\widetilde P$ they are. We can think of this attack relation in $\widetilde P$ as being ``cancelled out'' by the decoration of $a$'s step. So all units of dinv are compensated exactly by the negative contributions associated to valleys.
\end{proof}
We can view $\widetilde{P}$ as a ``base point'' of the cutting cycle of $P$, and then introduce a cyclic order on the cutting cycle $\CC(P)$. Be mindful that we might have $\sched(P)=1^n$ but $\sched(\tilde P)\neq 1^n$. 

\begin{definition}\label{def:order-hsteps}
    Take $P\in\LSQ(n)^{\bullet k}$, and $\widetilde{P}$ its canonical representative. We define an order on the horizontal steps of $\widetilde{P}$ that are ``not part of a decorated valley'' (that is not immediately followed by a decorated valley) by traversing each diagonal from right to left, from lowest to highest diagonal. More precisely, if $\ell_i$ is the $i$th horizontal step in this order, we define $Q_i \coloneqq \psi_{\ell_i}(\widetilde{P})$ and $Q_0 \coloneqq \widetilde{P}$.
\end{definition}
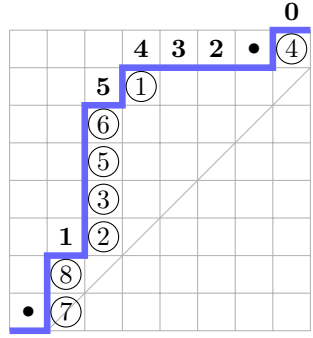
\begin{figure}
    \centering
    \begin{tikzpicture}[scale=.5]
        \draw[step=1.0, gray!60, thin] (0,0) grid (8,8);

        \draw[gray!60, thin] (1,0) -- (8,7);
        \node at (7.5,8.5) {\textbf{0}};
        \node at (1.5,2.5) {\textbf{1}};
        \node at (5.5,7.5) {\textbf{2}};
        \node at (4.5,7.5) {\textbf{3}};
        \node at (3.5,7.5) {\textbf{4}};
        \node at (2.5,6.5) {\textbf{5}};

        \draw[blue!60, line width=2.5pt] (0,0) -- (1,0) -- (1,1) -- (1,2) -- (2,2) -- (2,3) -- (2,4) -- (2,5) -- (2,6) -- (3,6) -- (3,7) -- (4,7) -- (5,7) -- (6,7) -- (7,7) -- (7,8) -- (8,8);

        \node at (1.5,0.5) {$7$};
        \draw (1.5,0.5) circle (.4cm);
        \node at (1.5,1.5) {$8$};
        \draw (1.5,1.5) circle (.4cm);
        \node at (2.5,2.5) {$2$};
        \draw (2.5,2.5) circle (.4cm);
        \node at (2.5,3.5) {$3$};
        \draw (2.5,3.5) circle (.4cm);
        \node at (2.5,4.5) {$5$};
        \draw (2.5,4.5) circle (.4cm);
        \node at (2.5,5.5) {$6$};
        \draw (2.5,5.5) circle (.4cm);
        \node at (3.5,6.5) {$1$};
        \draw (3.5,6.5) circle (.4cm);
        \node at (7.5,7.5) {$4$};
        \draw (7.5,7.5) circle (.4cm);
        \node at (0--1-0.5,0+0.5) {$\bullet$};
        \node at (7-0-0.5,7+0.5) {$\bullet$};
    \end{tikzpicture}
    \caption{Illustration of Definition~\ref{def:order-hsteps}}\label{fig:order-hsteps}
\end{figure}

See Figure~\ref{fig:order-hsteps} to see the ordering of the horizontal steps of $\tilde P$ of Example~\ref{ex:Ptilde}. See also Figure~\ref{fig:order-hsteps-big} for a larger (schematic) example.

Our next goal is to show that this choice of ordering of $\CC(\widetilde{P})$ is a well-chosen one with respect to the dinv. Let us illustrate what we mean using our running example.
\begin{example}
    Using Figure~\ref{fig:order-hsteps}, we can see that in Figure~\ref{fig:cutting-cycle-1n} we have 
    \begin{center}
        \begin{tabular}{c|cccccc}
            &$Q_0$&$Q_1$&$Q_2$&$Q_3$&$Q_4$&$Q_5$\\\hline
            Figure & \ref{fig:cc-d0}& \ref{fig:cc-d1}& \ref{fig:cc-P}& \ref{fig:cc-d3}& \ref{fig:cc-d4} &\ref{fig:cc-d5}\\
            $\dinv$ & 0&1&2&3&4&5.
        \end{tabular}
    \end{center}
    
\end{example}
\begin{figure}
    \centering
    \newcommand{\n}{21}
    \begin{tikzpicture}[scale = .4]
        \foreach \i in {0,1,...,\n}{
            \draw[gray,thin, densely dotted] (0,\i) --(\n-\i,\n);
        }
        \draw[blue!60, line width=2.5pt] 
        (0,0)-|(1,1)-|(2,2)-|(3,3)
        (3,3)|-(4,5)|-(5,9)|-(6,11)|-(7,15)
        (7,15)|- (8,16)|-(10,17)|-(12,18)|-(16,19)|-(20,20)|-(21,21) ;
        \node at (.5,.5) {$\bullet$};
        \node at (1.5,1.5) {$\bullet$};            
        \node at (.5,.5) {$\bullet$};            
        \node at (6.5,15.5) {$\bullet$};            
        \node at (7.5,16.5) {$\bullet$};           
        \node at (9.5,17.5) {$\bullet$};            
        \node at (11.5,18.5) {$\bullet$};
        \node at (15.5,19.5) {$\bullet$};            
        \node at (19.5,20.5) {$\bullet$}; 
        \node at (20.5,21.5) {0};
        \node at (2.5,2.5) {1};
        \node at (18.5,20.5) {2};
        \node at (3.5,5.5) {3};
        \node at (17.5,20.5) {4};
        \node at (16.5,20.5) {5};
        \node at (14.5,19.5) {6};
        \node at (4.5,9.5) {7};
        \node at (13.5,19.5) {8};
        \node at (5.5,11.5) {9};
        \node at (12.5,19.5) {10};
        \node at (10.5,18.5) {11};
        \node at (8.5,17.5) {12};
        \draw (1.5,.5) circle (.4cm);
        \draw (2.5,1.5) circle (.4cm);
        \draw (3.5,2.5) circle (.4cm);
        \draw (3.5,3.5) circle (.4cm);
        \draw (3.5,4.5) circle (.4cm);
        \draw (4.5,5.5) circle (.4cm);
        \draw (4.5,6.5) circle (.4cm);
        \draw (4.5,7.5) circle (.4cm);
        \draw (4.5,8.5) circle (.4cm);
        \draw (5.5,9.5) circle (.4cm);
        \draw (5.5,10.5) circle (.4cm);
        \draw (6.5,11.5) circle (.4cm);
        \draw (6.5,12.5) circle (.4cm);
        \draw (6.5,13.5) circle (.4cm);
        \draw (6.5,14.5) circle (.4cm);
        \draw (7.5,15.5) circle (.4cm);
        \draw (8.5,16.5) circle (.4cm);
        \draw (10.5,17.5) circle (.4cm);
        \draw (12.5,18.5) circle (.4cm);
        \draw (16.5,19.5) circle (.4cm);
        \draw (20.5,20.5) circle (.4cm);
    \end{tikzpicture}
    \caption{A (schematic) canonical representative of a schedule $1^n$ square path and the ordering of its steps as in Definition~\ref{def:order-hsteps}}\label{fig:order-hsteps-big}.
\end{figure}
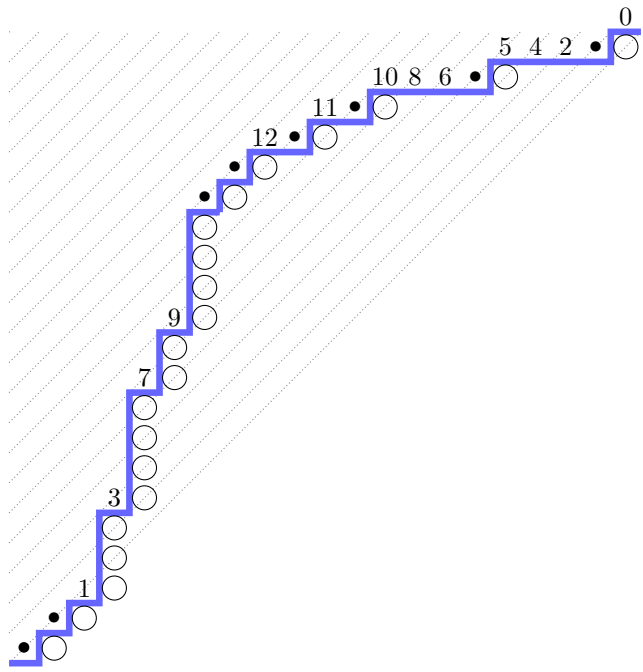
\begin{lemma}\label{CCDinvLemma}
    Take $P\in\stLSQ(n)^{\bullet k}$ such that $\sched(P) = 1^n$ and for any $0\leq i \leq n-k-1$ define $Q_i$ as in Definition~\ref{def:order-hsteps}. We have $\dinv(Q_i) = i$.
\end{lemma}
\begin{proof}
    By Lemma~\ref{lem:shape}, we know that each diagonal of $P$ (and hence of $Q_i$ for all $i$) contains at most two undecorated labels. Moreover, if a diagonal of $Q_0$ contains two undecorated labels, those labels must occur in consecutive columns (since the section of $Q_0$ containing the undecorated up steps contains no two consecutive east steps). This means that if a diagonal of $Q_0$ contains two undecorated labels then the leftmost undecorated label will be immediately followed by an east step. This east step will also necessarily be the second horizontal step in its diagonal in the ordering from~\ref{def:order-hsteps}.
    
    By the proof of Lemma~\ref{PTildedinv0}, we know that the only attack relations in $Q_0$ are of the form $(a,\dot{b})$ with $a$ to the left of $b$ and $a<b$. In particular, this means that in each diagonal of $\widetilde{P}$ the undecorated elements are in decreasing order when read left-to-right. Moreover, notice that as sets the diagonals of $Q_i$ are the same as those of $Q_0$; however, $\mathrm{shift}(Q_i) \geq \mathrm{shift}(Q_{i+1})$, and the elements of the diagonals may be permuted as we move from $Q_i$ to $Q_{i+1}$.

    Since we know that $\dinv(Q_0) = 0$, it suffices to show that $\dinv(Q_i) = \dinv(Q_{i-1}) + 1$ for all $i$. By our discussion in the previous paragraph, we know that each diagonal of $Q_0$ (as a set) appears as the $0$-diagonal of $Q_j$ and (possibly) $Q_{j+1}$ for some $j$. We consider multiple cases. We invite the reader to follow along on the schematic example in Figure~\ref{fig:order-hsteps-big}, which we take to be $Q_0$.
    
    \newlength{\case}
    \settowidth{\case}{Case 2.}
    \newlength{\subcase}
    \settowidth{\subcase}{Case 2.2.}
    \begin{enumerate}[labelwidth=\case, leftmargin=!]
        \item [\textbf{Case 1.}] Suppose that $Q_i$ and $Q_{i-1}$ have the same shift. For example, in Figure~\ref{fig:order-hsteps-big}, this happens for $i=1,3,7,9$.
        
        This happens exactly when the $0$ diagonal of $Q_{i-1}$ has two undecorated labels $a < b$. By our discussion above, we know that in $Q_{i-1}$ $b$ must be immediately to the left of $a$. Moreover, we know that the horizontal step immediately following the up step labeled $b$ will be the last step of $Q_i$. If we break $Q_{i-1}$ at the end of the horizontal step immediately following $b$ and write $Q_{i-1} = R_1 R_2$, then as in our definition of the cutting and pasting operation we have that $Q_{i} = R_2 R_1$. Since $Q_{i-1}$ and $Q_i$ have the same shift, we then see that the only contributions to $\dinv(Q_i)$ that do not also contribute to $\dinv(Q_{i-1})$ must come from pairs $(x,y)$ with $x \in R_2$ and $y \in R_1$. Certainly, $(a,b)$ is one such pair. Since $a$ is to the left of $b$ in $Q_i$, the pair now contributes $+1$ to $\dinv(Q_i)$. On the other hand, we claim that no other such pair can exist. 
        
        Indeed, suppose $(x,y)$ is such an attack pair, then since $x$ is to the left of $y$ in $Q_i$ it must be undecorated. In $Q_0$ and hence in $Q_i$, $x$ must be in a diagonal strictly above the one containing $a$ and $b$; and $y\in R_1$ in a diagonal strictly below the one containing $a,b$. It follows that $x$ and $y$ are at least two diagonals apart (in $Q_0$ and so also in $Q_i$) and do not form an attack relation. 
        We may deduce that $\dinv(Q_i) = \dinv(Q_{i-1}) + 1$.
        
        \item[\textbf{Case 2.}] Suppose that $\shift(Q_i) = \shift(Q_{i-1}) + 1$.
        \begin{enumerate}
            \item [\textbf{Case 2.1.}] Suppose that $Q_{i-1}$ had exactly one undecorated label $a$ in its zero diagonal. For example, in Figure~\ref{fig:order-hsteps-big}, this happens for $i=5,6,11,12$.
 
            We claim that the only new contribution to $\dinv(Q_i)$ is the contribution of $a$ to tertiary $\dinv$. Indeed, in this case the horizontal step labeled $i-1$ must be to the right of the rightmost undecorated up step (i.e., it must belong to the third segment of $Q_0$). This means that the only labels that have ``wrapped around'' (i.e., those that appeared after the horizontal step labeled $i$ in $Q_{i-1}$) when moving from $Q_{i-1}$ to $Q_i$ are decorated and in $Q_{i-1}$'s zero diagonal. Suppose $y$ is such a decorated label. In $Q_0$, we know that there is exactly one undecorated label $x$ to the left of $y$ such that $(x,y)$ is in attack position (which might be $a$). By construction, in $Q_{i-1}$ there is still exactly one undecorated $x$ to the left of $y$ with $(x,y)$ in attack position (the relative positions of all undecorated labels in $Q_{i-1}$ are the same as in $Q_0$). When we pass to $Q_i$, $y$ is to the left of $x$, so the contribution to $\dinv$ of the pair $(x,y)$ is unchanged --- $(x,y)$ is no longer in attack position, but $y$ is now in a negative diagonal. However, since $y$ is decorated, the net contribution to $\dinv$ in $Q_i$ is $0$ (as it was in $Q_{i-1}$). On the other hand, when we pass from $Q_{i-1}$ to $Q_i$, the unique undecorated label $a$ on the $0$ diagonal of $Q_{i-1}$ now lies on a negative diagonal and hence contributes $+1$ to the tertiary dinv of $Q_{i}$.
            \item[\textbf{Case 2.2.}] Finally, suppose that $Q_{i-1}$ has exactly two undecorated labels $a < b$ on its $0$ diagonal. For example, in Figure~\ref{fig:order-hsteps-big}, this happens for $i=2,4,8,10$.
            
            As we argued above, in $Q_{i-1}$ we will have $a$ to the left of $b$ since $b$ will be immediately followed by the final horizontal step of $Q_{i-1}$, and that the primary dinv of the pair $(a,b)$ is the only contribution to $\dinv(Q_{i-1})$ that did not also contribute to $\dinv(Q_{i-2})$. When we pass from $Q_{i-1}$ to $Q_{i}$ we can instead imagine that we skipped $Q_{i-1}$ and passed directly from $Q_{i-2}$ to $Q_{i}$. As we just argued, the only net contribution to $\dinv(Q_i)$ of this move comes from $a$ and $b$ moving from the $0$ diagonal of $Q_{i-2}$ to the $-1$ diagonal of $Q_i$. Thus, $\dinv(Q_i) = \dinv(Q_{i-2})+2 = \dinv(Q_{i-1})+1$, as desired.
        \end{enumerate} 
    \end{enumerate}
\end{proof}

\begin{lemma}\label{lem:CCCancel}
    If $P\in \stLSQ(n)^{\bullet k}$ such that $\sched(P)=1^n$ then \[\sum_{\substack{Q\in \CC(P)\\ \sched(Q) = 1^n}} (-1)^{\dinv(Q)}
        = \begin{cases}
            1 & \text{if $n-k$ is odd}  \\
            0 & \text{if $n-k$ is even}
        \end{cases}.\]
\end{lemma}

\begin{proof}
    By Observation~\ref{obs:cc}, Lemma~\ref{PTildedinv0}, and Lemma~\ref{CCDinvLemma}, we know that:
    \begin{itemize}
        \item $\CC(P)$ contains $n-k$ elements, and
        \item there exists an ordering $Q_0,\dots,Q_{n-k-1}$ of the elements of $\CC(P)$ such that $\dinv(Q_i) = i$.
    \end{itemize}
    In particular, this tells us that $\displaystyle \sum_{Q\in\CC(P)} q^{\dinv(Q)} = [n-k]_q$. Moreover, by the proof of Lemma~\ref{CCDinvLemma} and Proposition~\ref{prop:sched1-in-cc}, the $Q_i$ in $\CC(P)$ whose schedule words are not $1^n$ are exactly those $Q_i$ with two undecorated elements on the $0$-diagonal. Since all such $Q_i$ occur in consecutive pairs and $\dinv(Q_i) = i$, we have that
    $$\left.[n-k]_q\right|_{q=-1} =  \sum_{Q \in \CC(P)}(-1)^{\dinv(Q)} = \sum_{\substack{Q \in \CC(P)\\ \sched(Q) = 1^n}}(-1)^{\dinv(Q)},$$
    which yields the desired result.
\end{proof}

\begin{proof}[Proof of Theorem~\ref{thm:main1.2}]
Recall that by Equation~(\ref{eq:sched1sum}) 
$$S_{n,k} = \sum_{\substack{P \in \stLSQ(n)^{\bullet k}\\ \sched(P) = 1^n}}(-1)^{\dinv(P)}t^{\area(P)}.$$

Now, since all $Q$ in $\CC(P)$ have the same area, we can rewrite this sum as
\begin{align*}
    \sum_{\substack{P \in \stLSQ(n)^{\bullet k}\\ \sched(P) = 1^n}}(-1)^{\dinv(P)}t^{\area(P)} &= \sum_{C \in \mathcal{S}_{n,k}}t^{\area(C)}\sum_{\substack{Q \in \CC(C)\\ \sched(Q) = 1^n}}(-1)^{\dinv(Q)}\\
    &= \begin{cases}\displaystyle\sum_{C \in \mathcal{S}_{n,k}}t^{\area(C)}, & n-k \text{ odd}\\
    0, & n-k \text{ even}\end{cases},
\end{align*}
with the last equality following from Lemma~\ref{lem:CCCancel}.
\end{proof}

\begin{proposition}
    Theorem~\ref{thm:main1.2} is equivalent to Theorem~\ref{thm:main1}.
\end{proposition}
\begin{proof}
    The elements of $ \mathcal S_{n,k}$ are in bijection with those of $\ADR_{n,k}$.
    Indeed, given $C\in \mathcal S_{n,k}$, we map it to $\dw(C)$, which is clearly in $\ADR_{n,k}$ since at least one path in $C$ has schedule $1^n$. Conversely, given $\sigma\in \ADR_{n,k}$, we know that there is at least one shift $s$ for which $\sched(\sigma,s) = 1^n$. Given any such $s$ we know by Corollary~\ref{cor:schedule} that there is a unique $P \in \stLSQ(n)^{\bullet k}$ with $\sdw(P) = (\sigma,s)$. Suppose that $s_1$ and $s_2$ are shifts with $s_1 < s_2$ such that $\sched(\sigma,s_1) = \sched(\sigma,s_2) = 1^n$, and let $P_1$ and $P_2$ be the (unique) paths in $\stLSQ(n)^{\bullet k}$ with $\sdw(P_i) = (\sigma, s_i)$. We claim that by Corollary~\ref{cor:shape-canonical}, we must have $\tilde P_1 = \tilde P_2$. 
    
    Indeed, both paths start with a sequence of decorated labels in the $-1$-diagonal (which is nonempty if and only if the first decreasing run of $\sigma$ has no undecorated letters). There is no choice in the ordering of these valleys since they must be contractible. 
    
    Next, starting in the zero diagonal, come the undecorated labels, each in their own diagonal which may be deduced from $\sigma$. Each diagonal has either one or two undecorated labels; and if there are two they must be in a non-attack position.
    
    Finally, we place the remaining decorated labels in the appropriate diagonal in the last stretch. Again, since the valleys must be contractible there is no choice in the relative order of decorated valleys in the same diagonal. 

    Thus, since both paths have the same canonical representative, they must be part of the same cutting cycle.

\end{proof}

\section{From Dyck paths to square paths}\label{sec:dycktosquare}
The enumerators and combinatorics of $S_{n,k}$ and $D_{n,k}$ are closely connected. We explore these relationships in this section.

\subsection{A bijection}
In Table~\ref{tab:polvalues}, one may notice the following pattern.
\begin{proposition}\label{prop:sum-k-nfactorial}
    For all $n$, we have \[\sum_{k=0}^{n-1} S_{n,k} = \sum_{k=0}^{n-1} D_{n,k} = [n]_t!\]
\end{proposition}
It is established in the twin theorems \ref{thm:main1} and \ref{thm:main1.2}, that $S_{n,k}=0$ when $n-k$ is even, so all the contributions to the sum $\sum_k S_{n,k}$ must be contained in the terms where $n-k$ is odd.
The following lemma, combined with Theorem~\ref{thm:main1} and Theorem~\ref{thm:Dnk-enumeration}, gives an explicit combinatorial explanation for Proposition~\ref{prop:sum-k-nfactorial}.
It is a generalization of \cite{CorteelJosuatVergesVandenWyngaerd2023}*{Lemma~3.17}. In other words, it gives an explicit bijection of the combinatorics at $q=-1$ of the Hilbert series of the valley Delta conjecture and the valley Delta square conjecture. 
\begin{proposition}\label{prop:decorating-algos}
    For each permutation $\sigma\in\mathfrak S_{n}$, there exists
    \begin{enumerate}
        \item exactly one Dyck ADR\footnote{We recall that ADR stands for \emph{alternating Dyck representative}, see Definition~\ref{def:ADR}.} with underlying permutation $\sigma$;
        \item exactly one ADR with underlying permutation $\sigma$, such that its number of undecorated letters is odd.
    \end{enumerate}
\end{proposition}
\begin{proof}
    In this proof we will make extensive use of the left/right maximal cyclic run formulation of schedule numbers. In particular, we recall that LMCR and RMCR stand for \emph{left maximal cyclic run} and \emph{right maximal cyclic run}, respectively, as defined Definition~\ref{def:cyclic-run}.

    We will describe two slightly different algorithms, that for each permutation $\sigma\in\mathfrak S_n$, outputs a way of decorating $\sigma$ such that we obtain either a Dyck ADR or an ADR with an odd number of undecorated letters. We will call these algorithms the \emph{Dyck decorating algorithm} and the \emph{parity decorating algorithm}, respectively. The former was already described in \cite{CorteelJosuatVergesVandenWyngaerd2023}*{Lemma~3.17}, we re-state it here, for the sake of completeness.

    Let us make a few preliminary observations.
    \begin{itemize}
        \item If the last letter $\sigma_n$ of a permutation is decorated, its schedule number is zero.
        \item The $LMCR(\sigma_j)$ contains at least two letters for all $j>1$.
        \item  If $i<j$ and $\sigma_i\in LMCR(\sigma_j)$, if and only if $\sigma_j\in RMCR(\sigma_i)$.
    \end{itemize}

    Let us now describe the algorithms, whose first steps are in common (the reader is invited to follow along with Example~\ref{ex:decorating-algos}):

    \begin{enumerate}
        \item Set $j=n$.
        \item If $j>1$, consider $LMCR(\sigma_j)=\sigma_i\cdots \sigma_j$ for some $i<j$. Decorate all these letters except $\sigma_i$ and $\sigma_j$. If $j=1$, go to step 4.
        \item Set $j=i$ and start again at step 2.
    \end{enumerate}
    Only the last step differs for the two algorithms. 
    \begin{enumerate}
        \item[4.] (Dyck decorating algorithm) If the first decreasing run contains two undecorated letters, decorate $\sigma_1$.
        \item[4.] (Parity decorating algorithm) If the number of undecorated letters in $\sigma$ is even, decorate $\sigma_1$.
    \end{enumerate}
    In both algorithms, we end up with a decorated permutation such that there exists at least one run with exactly one undecorated letter. Indeed, in the Dyck decorating algorithm, the first run is such a run. In the parity decorating algorithm, since the total number of undecorated letters is odd and each run has at most two undecorated letters there must be a run with exactly one. We will show that, when choosing \textbf{any} shift $s$ such that a run with exactly one decorated letter becomes the zero run, all its shifted schedules are equal to $1$. Thus, we will have shown that the output of both algorithms is an ADR, and the output of the Dyck decorating algorithm is a Dyck ADR.

    First, notice that the schedule of the unique undecorated letter in the zero run is $1$.

    Next, it is clear from the definition of the schedule numbers that by construction, undecorated letters in positive runs have schedule number $1$.

    Let $\sigma_i$ be a decorated letter, or an undecorated letter in a negative run. We show that its shifted schedule number must be $1$.
    \begin{enumerate}
        \item If $w_{\sigma,s}(\sigma_i)\geq 2$ then there exists $j,k$ with $i<j<k$ such that $\sigma_j,\sigma_k\in RMCR(\sigma_i)$ and are undecorated. But we can see that $LMCR(\sigma_k)$ must then contain at least $\sigma_i\cdots \sigma_j\cdots \sigma_k$. This is contrary to the decorating algorithms since by construction the LMCR the undecorated letter $\sigma_k$ does not extend beyond the next undecorated letter to its left.
        \item If $w_{\sigma,s}(\sigma_i)=0$, take $j$ to be minimal such that $i<j$ and $\sigma_j$ is undecorated. Then if $LMCR(\sigma_j)=\sigma_k\cdots \sigma_j$, we must have $k\leq i$, because all letters between $\sigma_i$ and $\sigma_j$ are decorated. But then $\sigma_j\in RMCR(\sigma_i)$ and so its schedule cannot be $0$.
    \end{enumerate}

    It is left to show uniqueness. Let $\sigma$ be an ADR with an odd number of decorated letters or a Dyck ADR. We will show that the decorations on $\sigma$ must be exactly the ones determined by the appropriate decorating algorithm. Fix a shift $s$ for which $w_{\sigma,s}(\sigma_i)=1$ for all $i$.

    We must have that $\sigma_n$ is undecorated because otherwise its schedule is zero.

    Let $\sigma_j$ be an undecorated letter in a positive run, then $LMCR(\sigma_j)=\sigma_i\cdots \sigma_j$ contains exactly one undecorated letter, say $\sigma_k$ with $k<j$. If $i<k$ then $\sigma_i$ is decorated and $RMCR(\sigma_i)$ contains the undecorated $\sigma_k$ and $\sigma_j$ so $w_{\sigma,s}(\sigma_i)\geq 2$. So we must have $i=k$, as prescribed by the algorithm.

    Let $\sigma_i$ be an undecorated letter in a negative run and $j$ minimal such that $i<j$ and $\sigma_j$ is undecorated. We will show that $LMCR(\sigma_j)=\sigma_i\cdots \sigma_j$.
    \begin{itemize}
        \item If $\sigma_i\not\in LMCR(\sigma_j)$ then
              \begin{itemize}
                  \item either there is a decreasing run between $\sigma_i$'s and $\sigma_j$'s run. In this case there is a run without undecorated letters which is not the first run. To see that this is not possible, one can argue that in the corresponding square path, diagonals that are not the bottom one must contain a vertical step which is preceded by another vertical step, which may not be decorated;
                  \item or $\sigma_i$ is in the run preceding $\sigma_j$'s run and $\sigma_i>\sigma_j$. It follows that $RMCR(\sigma_i)=\sigma_i\cdots \sigma_k$ for some $k<j$ and so $w_{\sigma,s}(\sigma_i) = 0$.
              \end{itemize}
        \item If $\sigma_{i-1}\in LMCR(\sigma_j)$ then $\sigma_i, \sigma_j\in RMCR(\sigma_{i-1})$ and so $w_{\sigma_{i-1}} \geq 2$.
    \end{itemize}
    Thus, since there can be only one undecorated letter in the zero run, we have shown that if $i<j$ with $\sigma_i,\sigma_j$ undecorated and $\sigma_k$ decorated for all $i<k<j$, then $LMCR(\sigma_j)=\sigma_i\cdots\sigma_j$; and that $\sigma_n$ is undecorated. So the decorations on $\sigma_2\cdots \sigma_n$ are as prescribed by the common steps 1-3 of both algorithms.
    Finally, given these forced decorations, the decoration condition on $\sigma_1$ exactly reduces to step 4 in the respective algorithms.
\end{proof}
Clearly, this result, in combination with Theorem~\ref{thm:main1}, yields the predicted condensed formula \eqref{eq:sum-Snkz} from the introduction. 

\begin{example}\label{ex:decorating-algos}
    Take $n=9$ and consider the permutation $852961743$. The left maximal cyclic runs to consider are $1743$, $961$, $8529$. Thus, after the common steps 1--3 of the decorating algorithms, we have $\sigma = 8\dec 5\dec 29\dec 61\dec 7\dec 43$.

    For the Dyck decorating algorithm, since the first decreasing run contains exactly one undecorated letter, we do not decorate $\sigma_1$ and keep $\sigma = 8\dec 5\dec 29\dec 61\dec 7\dec 43$.

    For the parity decorating algorithm, since there are an even number of undecorated letters, we must decorate the first letter and obtain $\sigma = \dec 8\dec 5\dec 29\dec 61\dec 7\dec 43$. One can see that for shift $1$ or $2$, this decorated permutation has schedules all $1$. 
\end{example}
Notice that in the proof of Proposition~\ref{prop:decorating-algos}, for a given $\sigma$, we might have chosen \textbf{any} shift $s$ such that the zero-run contains exactly one undecorated letter. 
\begin{corollary}\label{cor:anyshift}
    Let $\sigma\in \ADR_{n,k}$. Then the set of shifts $s$ for which the schedule numbers of $(\sigma,s)$ are all one is exactly the set of shifts such that the zero run contains exactly one undecorated letter. 
\end{corollary}
We can now easily establish a result we needed in the previous section. 
\begin{proof}[Proof of Proposition~\ref{prop:sched1-in-cc}]
    Take $P\in \LSQ(n)^{\bullet k}$ such that $\sched(P) = 1^n$ and $Q$ an element in its cutting cycle. We must have that $\sdw(P) = (\sigma,s_1)$ and $\sdw(Q) = (\sigma,s_2)$ for some decorated permutation $\sigma$ and integers $s_1,s_2$. Since $\sched(P) = 1^n$, $\sigma$ is an ADR. From Corollary~\ref{cor:anyshift} it follows that $Q$ has schedule $1^n$ if and only if it has one undecorated letter in the zero diagonal. 
\end{proof}

The following is deduced easily from Proposition~\ref{prop:decorating-algos}, by composing the bijections with $\mathfrak S_n$ there described.
\begin{corollary}
    There exists a bijection
    \begin{align*}
        \phi: \bigsqcup_{\substack{k\in [n-1] \\n-k\text{ is odd}}}\ADR_{n,k} \rightarrow \bigsqcup_{k\in[n-1]}\DADR_{n,k}.
    \end{align*}
    The underlying permutations of $\sigma$ and $\phi(\sigma)$ coincide, so $\revmaj(\sigma) = \revmaj(\phi(\sigma))$.
\end{corollary}
We can describe $\phi$ explicitly as follows. Take $\sigma$ and ADR with an odd number of decorations. If the first decreasing run of $\sigma$
\begin{itemize}
    \item contains no undecorated letters, remove the decoration from $\sigma_1$;
    \item contains exactly one undecorated letter, do nothing;
    \item contains two undecorated letters, decorate $\sigma_1$.
\end{itemize}

\begin{example}
    For $n=3$ the bijection is given by the following correspondence
    \begin{center}
        \begingroup
        \renewcommand{\arraystretch}{2}
        \begin{tabular}{c|cccccc}
            $\sigma$       & $123$ & $231$ & $\dec 1\dec 32$ & $312$       & $\dec 2 \dec 1 3$ & $\dec 3 \dec 2 1$  \\ \hline
            $\phi(\sigma)$ & $123$ & $231$ & $1\dec 32$      & $\dec{3}12$ & $2\dec{1}3$       & $\dec{3}\dec{2}1$.
        \end{tabular}
        \endgroup
    \end{center}
\end{example}
\subsection{A recursion}
In this subsection, we explore a recursive relationship between the $S_{n,k}$ and the $D_{n,k}$. The advantage of this recursive approach is that we can understand $S_{n,k}$ for a fixed $k$, and not merely when taking the sum over $k$, as in the previous subsection. The recursive generation is similar to the one described in \cite{CorteelJosuatVergesVandenWyngaerd2023}*{Section~5.1}.

\begin{theorem}\label{thm:dycktosquare}
    For all $n,k\in \N$, we have
    \[
        S_{n,k} = \begin{cases}
            [n]_t\left(D_{n-1,k}+D_{n-1,k-1}\right) & \text{if $n-k$ is odd}  \\
            0                                       & \text{if $n-k$ is even}
        \end{cases},
    \]where we use the convention that $D_{-1,k} = 0$ for all $k$.
\end{theorem}

\begin{proof}[Proof of Theorem~\ref{thm:dycktosquare}]
    We have already established, in the twin theorems \ref{thm:main1} and \ref{thm:main1.2}, that $S_{n,k}=0$ when $n-k$ is even. Let us thus suppose that $n-k$ is odd, in which case we have
    \[S_{n,k} = \sum_{\tau\in \ADR_{n,k}}t^{\revmaj(\tau)}.\]

    Take $\sigma \in \mathfrak \DADR_{n-1,k}\sqcup \DADR_{n-1,k-1}$, and $m\in \{1,\dots,n\}$. 
    We define $\delta_m(\sigma)\in\mathfrak S_n^{\bullet k}$ such that its $i$-th letter is
    \begin{itemize}
        \item $m$ if $i=1$;
        \item $m + \sigma_{i-1} \pmod{n}$  if $i>1$, where we take the representative modulo $n$ in $\{1,\dots, n\}$.   
    \end{itemize}
    The decorations on $\delta_m(\sigma)$ are such that if the $i$-th letter is decorated in $\sigma$, then the $(i+1)$-th letter is decorated in $\delta_m(\sigma)$. Furthermore, if $\sigma\in\DADR_{n-1,k-1}$, decorate $\sigma_1$. See Example~\ref{ex:cylinder-cut} and its Figure~\ref{fig:cylinder-cut} for a helpful visual interpretation of this map. 
    
    Let us show that for all $m\in \{1,\dots,n\}$, $\tau\coloneqq \delta_m(\sigma)\in \ADR_{n,k}$. We do this by showing that the decorating conditions of $\delta_m$ coincide with those described by the parity decorating algorithm (Proposition~\ref{prop:decorating-algos}).

    It is clear from the definition of a cyclic run that if $\sigma_i\cdots \sigma_j$ for $1\leq i\leq j\leq n$ is contained in a cyclic run, then so is $\tau_{i+1}\cdots \tau_{j+1}$. It follows that in the (parity) decorating algorithm for $\tau$, the decorations on $\tau_3\cdots\tau_n$ coincide with those on $\sigma_2\cdots \sigma_{n-1}$, obtained from the Dyck decorating algorithm (since both algorithms coincide on all but the first letter and the cyclic runs are the same). Let $j>1$ be minimal such that $\sigma_j$ is undecorated. Then $\sigma_1\cdots \sigma_j$ is a cyclic run, so in the (parity) decorating algorithm $\tau_2$ is decorated if and only if $\tau_1\tau_2\cdots \tau_{j+1}$ is a cyclic run. Knowing that $\tau_2\cdots \tau_{j+1}$ is a cyclic run, this is the case if and only if one of the following is true:
    \begin{align}
        m=\tau_1>\tau_2>\tau_{j+1}\label{eq:maxrun1}\\
        m=\tau_1>\tau_2\text{ and }\tau_1<\tau_{j+1}\label{eq:maxrun2}\\
        m=\tau_1<\tau_{j+1}\text{ and }\tau_2>\tau_{j+1}\label{eq:maxrun3}
    \end{align}
    We show that this condition (for $\tau_2$ to be decorated in the parity decorating algorithm of $\tau$) coincides with the condition for $\sigma_1$ to be decorated in the Dyck decorating algorithm. Indeed, $\sigma_1$ is decorated in the Dyck decorating algorithm if and only if the first decreasing run of $\sigma$ contained two undecorated letters after step 3 of the algorithm. This is the case if and only if $\sigma_j$ is in $\sigma$'s first decreasing run, i.e. if and only if \ $\sigma_1>\sigma_j$. Let us consider three cases. 
    \begin{enumerate}
        \item If $n-m<\sigma_j<\sigma_1$, then $\tau_2=\sigma_1+m-n$ and $\tau_{j+1} = \sigma_j+m-n$; and we have \eqref{eq:maxrun1}.
        \item If $\sigma_j\leq n-m < \sigma_1$, then $\tau_2 = \sigma_1+m-n$ and $\tau_{j+1} = \sigma_j+m$; and we have \eqref{eq:maxrun2}.
        \item If $\sigma_j<\sigma_1\leq n-m$, then $\tau_2 = \sigma_1+m$ and $\tau_{j+1} = \sigma_j+m$; and we have \eqref{eq:maxrun3}.
    \end{enumerate}
    It follows that the parity decorating algorithm on $\tau_2\dots\tau_{n}$ produces the same decorations as the Dyck decorating algorithm on $\sigma$ (shifted by one index). Finally, the parity decorating algorithm prescribes that $\tau_1$ is decorated if and only if the number of undecorated letters in $\tau_2\cdots \tau_n$ is even, i.e.\ if and only if $\sigma\in\DADR_{n-1,k-1}$ (since we supposed $n-k$ to be odd).

    Thus, we have shown that, if $\sigma \in \DADR_{n-1,k}\sqcup \DADR_{n-1,k-1}$, $\delta_m(\sigma)\in \ADR_{n,k}$, for all $m\in\{1,\dots,n\}$. Since $\sqcup_{k=0}^{n-1}\DADR_{n-1,k}$ consists of all elements of $\mathfrak S_{n-1}$ (with the appropriate decorations), it is clear that by applying $\delta_m$ to all its elements for all $m\in [n]$, we obtain all elements in $\mathfrak S_{n}$, with the appropriate decorations. In other words, this process generates all elements if $\sqcup_{k=0}^{n-1}\ADR_{n,k}$ (see Example~\ref{ex:recursion}). 

    It is left to show that we get the right $t$-contributions. For this, it suffices to notice that \[\revmaj(\delta_m(\sigma)) = \revmaj(\sigma) + n-m.\] The proof of this fact is the same as in \cite{CorteelJosuatVergesVandenWyngaerd2023}*{Proposition~5.9}. 
\end{proof}
\begin{example}\label{ex:cylinder-cut}
    We consider the following DADR (as obtained in Example~\ref{ex:decorating-algos}): \[\sigma = 8\,\dec 5\,\dec 2\,9\,\dec 6\,1\,\dec 7\,\dec 4\,3.\] Then \[\delta_4(\sigma) = 4\,2\,\dec 9\,\dec 6\,3\,\dec{10}\,5\,\dec 1\,\dec 8\,7.\]
    We can represent this operation visually, on a cylinder, as in Figure~\ref{fig:cylinder-cut}.
\end{example}
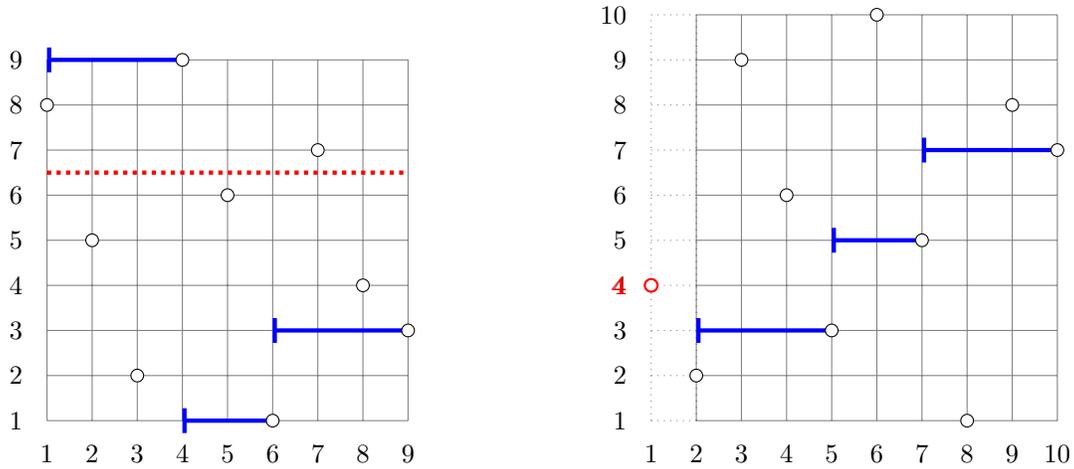
\begin{figure}[H]
    \centering
\begin{tikzpicture}[scale=.6]
    \draw[opacity = .5] (1,1) grid (9,9);
    \draw[ultra thick, blue,-{Bar[]}] (9,3) -- (6,3);
    \draw[ultra thick, blue,-{Bar[]}](6,1)--(4,1);
    \draw[ultra thick, blue,-{Bar[]}](4,9)--(1,9);
    \draw[ultra thick, dotted, red] (1,6.5) -- (9,6.5);
    \foreach \i/\j in {1/8,2/5,3/2,4/9,5/6,6/1,7/7,8/4,9/3}{ 
    \node[left=.2cm] at (1,\i) {\i};
    \node[below=.2cm] at (\i,1) {\i};
    \filldraw[fill = white] (\i,\j) circle (4pt);
    }        
\end{tikzpicture}
\hspace{2cm}
\begin{tikzpicture}[scale=.6]
    \draw[opacity = .5] (2,1) grid (10,10);
    \draw[opacity = .5,dotted] (1,1) grid (2,10);
    \draw[ultra thick, blue,-{Bar[]}] (10,7) -- (7,7);
    \draw[ultra thick, blue,-{Bar[]}](7,5)--(5,5);
    \draw[ultra thick, blue,-{Bar[]}](5,3)--(2,3);
    \foreach \i/\j in {1/4,2/2,3/9,4/6,5/3,6/10,7/5,8/1,9/8,10/7}{ 
    \node[left=.2cm] at (1,\i) {\i};
    \node[below=.2cm] at (\i,1) {\i};
    \filldraw[fill = white] (\i,\j) circle (4pt);
    }        
    \filldraw[fill=white,draw=red,thick] (1,4)circle(4pt) node[left=.2cm,red] {\textbf{4}};
\end{tikzpicture}
\caption{$\sigma$ and $\delta_4(\sigma)$ }
\label{fig:cylinder-cut}
\end{figure}

\begin{example}\label{ex:recursion}
    Take $n = 4$. In Table~\ref{tab:recursion}, we apply $\delta_m$ to the 6 elements of $\sqcup_{k=0}^{2}\DADR_{3,k}$, for $m=1,\dots, 4$, to obtain all 24 elements in $\sqcup_{k=0}^{3}\ADR_{4,k}$. 
    \begin{table}[H]
        \begin{center}
            \begingroup
            \renewcommand{\arraystretch}{1.5}
            \begin{tabular}[t]{cc}
            $\DADR_{3,0}$ & $\ADR_{4,1}$ \\ \hline\hline
            123           & $\dec{1}234$ \\
                          & $\dec{2}341$ \\
                          & $\dec{3}412$ \\             
                          & $\dec{4}123$ \\ \hline             
            231           & $\dec{1}342$ \\
            & $\dec{2}413$ \\
            & $\dec{3}124$ \\
                          & $\dec{4}231$ \\
        \end{tabular}
        \hspace{1cm}
        \begin{tabular}[t]{cc}
            $\DADR_{3,1}$ & $\ADR_{4,1}$ \\ \hline\hline
            $1\dec{3}2$   & $12\dec{4}3$ \\
            & $23\dec{1}4$ \\
            & $34\dec{2}1$ \\             
            & $41\dec{3}2$ \\ \hline             
            $\dec{3}21$   & $1\dec{4}32$ \\
            & $2\dec{1}43$ \\
            & $3\dec{2}14$ \\
            & $4\dec{3}21$ \\ \hline
            $2\dec{1}3$   & $13\dec{2}4$ \\
            & $24\dec{3}1$ \\
            & $31\dec{4}2$ \\
            & $42\dec{1}3$ \\
        \end{tabular}
        \hspace{1cm}
        \begin{tabular}[t]{cc}
            $\DADR_{3,2}$ & $\ADR_{4,3}$\\
            \hline\hline
            $\dec 3\dec 2 1$ & $\dec 1\dec 4\dec 32$\\
            & $\dec 2 \dec 1 \dec 43$\\
            & $\dec 3 \dec 2 \dec 1 4$\\ 
            & $\dec 4\dec 3\dec 21$\\ 
        \end{tabular}
        \endgroup
    \end{center}
    \caption{Generate $\sqcup_{k=0}^{3}\ADR_{4,k}$ from $\sqcup_{k=0}^{2}\DADR_{3,k}$}\label{tab:recursion}
    \end{table}
\end{example}


\section{Concluding remarks}
\begin{itemize}
    \item In \cite{IraciVandenWyngaerd2021}, the authors formulate a slightly ``modified'' version of the valley Delta square conjecture, using the novel Theta operators from \cite{DAdderioIraciVandenWyngaerd2021}: it is a formula for $\Theta_{e_{k}}\nabla \omega(p_{n-k})$ in terms of the subset of labelled decorated square paths that have at least one non-decorated vertical step on its lowest diagonal. The symmetric functions of both versions are linked by the following formula:
    \[\frac{[n]_t}{[n-k]_t}\Theta_{e_k}\nabla\omega(p_{n-k}) = \frac{[n-k]_q}{[n]_q}\Delta_{e_{n-k}}\omega(p_n). \]
    It seems that, evaluated at $q=-1$, the Hilbert series of $\Theta_{e_{k}}\nabla \omega(p_{n-k})$ is also always $t$-positive and equal $0$ whenever $n-k$ is even. A natural question would be to adapt the arguments of this paper to the ``modified'' setting. At first view however, the combinatorics seems to behave less nicely. In particular, the sum over $k$ does not yield $[n]_t!$, as before. 
    \item Inspired by Theorem~\ref{thm:dycktosquare}, we have the following conjectural symmetric function formula.
    \begin{conjecture} For all positive integers $n$ and partitions $\lambda$ of $n-1$, we have
        \begin{align*}
            &\left.\left\langle\frac{[n-k]_q}{[n]_q}\Delta_{e_{n-k}}\omega(p_{n}),h_{(\lambda,1)}\right\rangle \right|_{q=-1}\\
            &\hspace{1cm}=\begin{cases}
                [n]_t\left.\left(\left\langle \Delta'_{e_{(n-1)-k-1}}e_{n-1}+ \Delta'_{e_{(n-1)-(k-1)-1}}e_{n-1},h_{\lambda} \right\rangle\right)\right|_{q=-1}& \text{if $n-k$ is odd}\\
                0 &\text{if $n-k$ is even.}
            \end{cases}
        \end{align*}
    \end{conjecture}
    \item As mentioned in \cite{CorteelJosuatVergesVandenWyngaerd2023}, computational evidence suggests that the evaluation at $q =-1$ yields $t$-positive results for many other polynomials related to the shuffle theorem and Delta conjectures. We refer to the aforementioned paper for more details. 
\end{itemize}
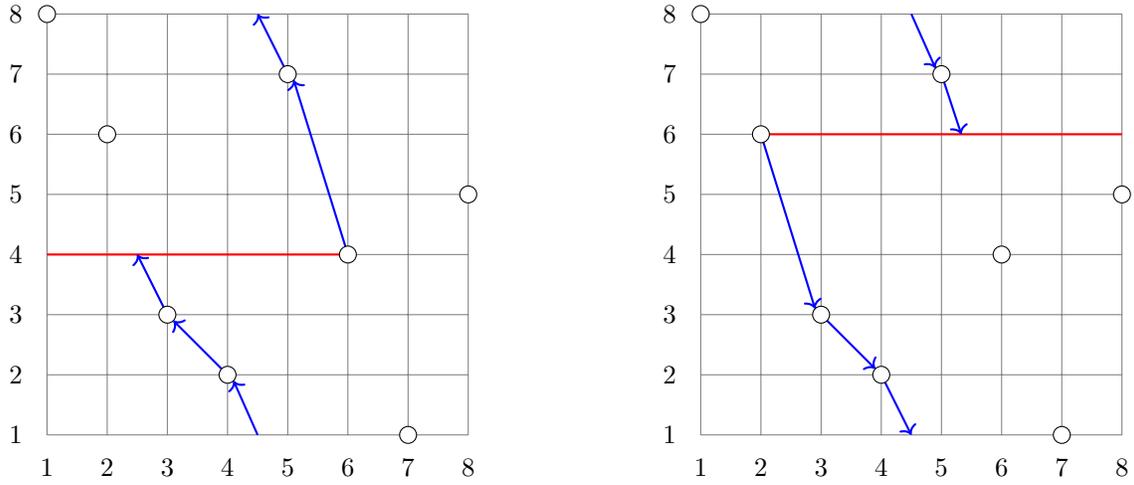
\begin{figure}[H]
    \centering
    \begin{tikzpicture}[scale=.8]
        \draw[opacity = .5] (1,1) grid (8,8);
        \foreach \i in {1,...,8}{
                \node[left=.2cm] at (1,\i) {\i};
                \node[below=.2cm] at (\i,1) {\i};
            }
        \draw[thick,red](6,4) -- (1,4);
        \draw[->,thick,blue] (6,4)--(5.1, 6.9);
        \draw[->,thick,blue] (5,7)--(4.5, 8);
        \draw[->,thick,blue] (4.5,1)--(4.1, 1.9);
        \draw[->,thick,blue] (4,2)--(3.1, 2.9);
        \draw[->,thick,blue] (3,3)--(2.5,4);
        \draw[fill= white] (1,8) circle (4pt);
        \draw[fill= white] (2,6) circle (4pt);
        \draw[fill = white] (3,3) circle (4pt);
        \draw[fill= white] (4,2) circle (4pt);
        \draw[fill = white] (5,7) circle (4pt);
        \draw[fill = white] (6,4) circle (4pt);
        \draw[fill= white] (7,1) circle (4pt);
        \draw[fill = white] (8,5) circle (4pt);
    \end{tikzpicture}
    \hfill
    \begin{tikzpicture}[scale=.8]
        \draw[opacity = .5] (1,1) grid (8,8);
        \foreach \i in {1,...,8}{
                \node[left=.2cm] at (1,\i) {\i};
                \node[below=.2cm] at (\i,1) {\i};
            }
        \draw[thick,red](2,6) -- (8,6);
        \draw[->, thick, blue] (2,6) -- (2.9,3.1);
        \draw[->, thick, blue] (2.9,3.1)--(3.9,2.1);
        \draw[->, thick, blue] (4,2)--(4.5,1);
        \draw[->, thick, blue] (4.5,8)--(4.9,7.1);
        \draw[->, thick, blue] (5,7)--(5.33,6);

        \draw[fill= white] (1,8) circle (4pt);
        \draw[fill= white] (2,6) circle (4pt);
        \draw[fill = white] (3,3) circle (4pt);
        \draw[fill= white] (4,2) circle (4pt);
        \draw[fill = white] (5,7) circle (4pt);
        \draw[fill = white] (6,4) circle (4pt);
        \draw[fill= white] (7,1) circle (4pt);
        \draw[fill = white] (8,5) circle (4pt);
    \end{tikzpicture}
    \caption[]{Illustration of left/right maximal cyclic runs of the permutation $86327415$.}\label{fig:cyclic-runs}
\end{figure}
\afterpage{%
\thispagestyle{empty}
\newgeometry{a4paper, headsep=0pt, footskip=0pt, inner=0pt, outer=0pt, bottom=10pt, top=10pt}
\begin{figure}[H]
    \centering
    \includegraphics[scale=.42]{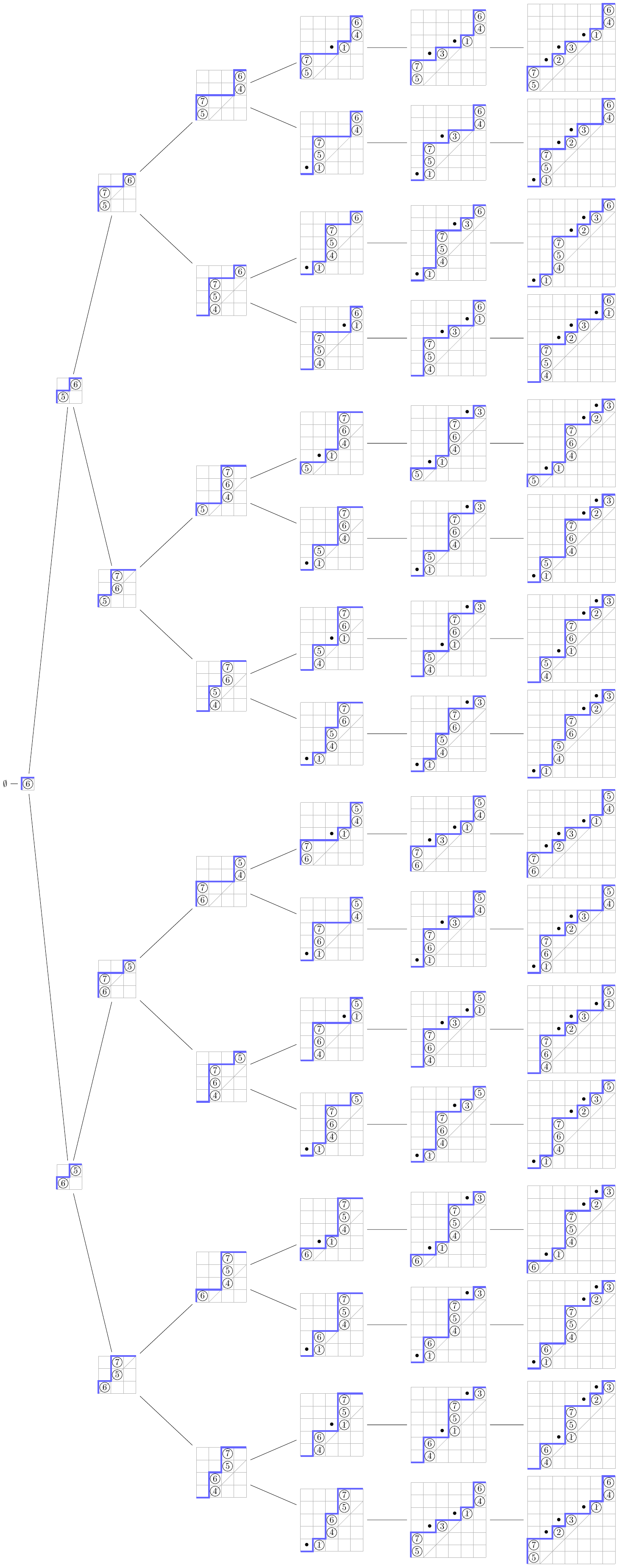}
    \caption{Illustration of the proof of Theorem~\ref{thm:sched-formula}}\label{fig:schedule-tree}
\end{figure}
\clearpage
\restoregeometry
}
\newpage

\bibliography{biblio}
\end{document}